\documentstyle[12pt]{article}

\def\baselinestretch{1} \topmargin -12pt \headsep 0pt \textwidth 470pt
\oddsidemargin 0pt \def\medskipamount{20pt} \def\arraycolsep{2pt}

\input prepictex
\input pictex
\input postpictex

\newcommand{\bit}[1]{\section{#1} \setcounter{equation}{0}}
\renewcommand{\theequation}{\thesection .\arabic{equation}}
\newcommand{\re}[1]{{\bf (\ref{#1})}}

\newenvironment{thm}{ \addvspace{\medskipamount} \refstepcounter{equation}
\noindent {\bf (\theequation)} \begin{em}}{\end{em} \par
\addvspace{\medskipamount} }

\newenvironment{lemma}{ \addvspace{\medskipamount} \refstepcounter{equation}
\noindent {\bf (\theequation)} \begin{em}}{\end{em} \par
\addvspace{\medskipamount} }

\newenvironment{propn}{ \addvspace{\medskipamount} \refstepcounter{equation}
\noindent {\bf (\theequation)} \begin{em}}{\end{em} \par
\addvspace{\medskipamount} }

\newenvironment{cor}{ \addvspace{\medskipamount} \refstepcounter{equation}
\noindent {\bf (\theequation)} \begin{em}}{\end{em} \par
\addvspace{\medskipamount} }

\newenvironment{defn}{ \addvspace{\medskipamount} \refstepcounter{equation}
\noindent {\bf (\theequation)} \/{\em Definition}.  }{ \par
\addvspace{\medskipamount} }

\newenvironment{remk}{ \addvspace{\medskipamount} \refstepcounter{equation}
\noindent {\bf (\theequation)} \/{\em Remark}.  }{ \par
\addvspace{\medskipamount} }

\newenvironment{manynotop}{ \addvspace{\smallskipamount}
\refstepcounter{equation}
\noindent }{ \par \addvspace{\medskipamount} }

\newfont{\bbb}{msbm10 scaled\magstep1}
\newfont{\blb}{msbm7 scaled\magstep1}

\begin{document}
\catcode`\@=\active
\catcode`\@=11


\def\@eqnnum{\hbox to .01pt{}\rlap{\bf \hskip -\displaywidth(\theequation)}}


\newcommand{\Bdal}[1]{{\phantom{}\Big\downarrow
\makebox[0pt][r]{$\scriptstyle #1\;\;\;$}}}

\newcommand{\Bdar}[1]{{\phantom{}\Big\downarrow
\raisebox{.4ex}{$\scriptstyle #1$}}}

\newcommand{\eg}{\mbox{$X_{g}$}}
\newcommand{\ega}{\mbox{$X$}}

\newcommand{\sltc}{\mbox{${\rm SL}(2,{\C})$}}
\newcommand{\sot}{\mbox{${\rm SO}(3)$}}
\newcommand{\psltc}{\mbox{${\rm PSL}(2,{\C})$}}
\newcommand{\cala}{\mbox{${\cal A}$}}
\newcommand{\calg}{\mbox{${\cal G}$}}

\newcommand{\hata}{\mbox{$\hat{\cal A}$}}
\newcommand{\hatf}{\mbox{$\hat{f}$}}
\newcommand{\hatfs}{\mbox{$\hat{f}^{*}$}}
\newcommand{\hatg}{\mbox{$\hat{\cal G}$}}
\newcommand{\hatga}{\mbox{$\hat{\cal G}_{a}$}}

\newcommand{\hatl}{\mbox{$\hat{\cal L}$}}
\newcommand{\call}{\mbox{${\cal L}$}}
\newcommand{\callk}{\mbox{${\cal L}^{k}$}}
\newcommand{\hatlk}{\mbox{$\hat{\cal L}^{k}$}}
\newcommand{\hatlka}{\mbox{$\hat{\cal L}^{k}_{a}$}}

\newcommand{\calng}{\mbox{$N$}}
\newcommand{\calo}{\mbox{${\cal O}$}}
\newcommand{\calu}{\mbox{${\cal U}$}}
\newcommand{\bc}{\mbox{${\C}$}}
\newcommand{\al}{\mbox{${\alpha}$}}

\newcommand{\ga}{\mbox{${\gamma}$}}
\newcommand{\gi}{\mbox{${\gamma_{i}}$}}
\newcommand{\la}{\lambda}
\newcommand{\si}{\mbox{${\sigma}$}}

\newcommand{\bu}{\mbox{${\bf U}$}}
\newcommand{\Zuk}{\mbox{$Z_{k}$}}
\newcommand{\zuk}{\mbox{$z_{k}$}}
\newcommand{\Ztk}{\mbox{$\hat{Z}_{k}$}}
\newcommand{\ztk}{\mbox{$\hat{z}_{k}$}}

\newcommand{\dgk}{\mbox{$D(g,k)$}}


\newcommand{\C}{\mbox{\bbb C}}
\newcommand{\Cs}{\mbox{\blb C}}
\newcommand{\Pj}{\mbox{\bbb P}}
\newcommand{\Ps}{\mbox{\blb P}}
\newcommand{\R}{\mbox{\bbb R}}
\newcommand{\Rs}{\mbox{\blb R}}
\newcommand{\Q}{\mbox{\bbb Q}}
\newcommand{\Qs}{\mbox{\blb Q}}
\newcommand{\Z}{\mbox{\bbb Z}}
\newcommand{\Zs}{\mbox{\blb Z}}
\newcommand{\N}{\mbox{\bbb N}}
\newcommand{\Ns}{\mbox{\blb N}}
\newcommand{\Coeff}{\mathop{{\rm Coeff}}}
\newcommand{\End}{\mathop{{\rm End}}\nolimits}
\newcommand{\Gr}{\mathop{{\rm Gr}}\nolimits}
\newcommand{\Hom}{\mathop{{\rm Hom}}\nolimits}
\newcommand{\La}{\Lambda}
\newcommand{\Pic}{\mathop{{\rm Pic}}\nolimits}
\newcommand{\Quot}{\mathop{{\rm Quot}}\nolimits}
\newcommand{\Res}{\mathop{{\rm Res}}}
\newcommand{\Spec}{\mathop{{\rm Spec}}\nolimits}
\newcommand{\SPEC}{\mathop{{{\bf Spec}}}\nolimits}
\newcommand{\beqas}{\begin{eqnarray*}}
\newcommand{\beqa}{\begin{eqnarray}}
\newcommand{\be}{\mbox{$\bf E$}}
\newcommand{\beq}{\begin{equation}}
\newcommand{\bet}{\mbox{$\tilde{\bf E}$}}
\newcommand{\bexl}{\be_{\bx}^{\bl}}
\newcommand{\bexlm}{\be_{\bx}^{\bl,\bm}}
\newcommand{\bira}{\leftrightarrow}
\newcommand{\bl}{{\bf \lll}}
\newcommand{\bla}{{\bf \Lambda}}
\newcommand{\bm}{{\bf \mmm}}
\newcommand{\bp}{{\bf \Phi}}
\newcommand{\bx}{{\bf q}}
\newcommand{\bxt}{\tilde{\bx}}
\newcommand{\ch}{\mathop{{\rm ch}}\nolimits}
\newcommand{\ci}{{\cal I}}
\newcommand{\clok}{\mbox{${\cal L}_{0}^{k}$}}
\newcommand{\clo}{\mbox{${\cal L}_{0}$}}
\newcommand{\cl}{\mbox{$\cal L$}}
\newcommand{\co}{{\cal O}}
\newcommand{\codim}{\mathop{{\rm codim}}\nolimits}
\newcommand{\coker}{\mathop{{\rm coker}}\nolimits}
\newcommand{\conga}{=}
\newcommand{\cu}{\mbox{$\cal U$}}
\newcommand{\cv}{\mbox{${\cal W}^{\bla}$}}
\newcommand{\cw}{\mbox{$\cal W$}}
\newcommand{\cxt}{\mbox{$\tilde{\cal X}$}}
\newcommand{\cx}{\mbox{$\cal X$}}
\newcommand{\cy}{\mbox{$\cal Y$}}
\newcommand{\diag}{\mathop{{\rm diag}}\nolimits}
\newcommand{\eeqas}{\end{eqnarray*}}
\newcommand{\eeqa}{\end{eqnarray}}
\newcommand{\eeq}{\end{equation}}
\newcommand{\eiiit}{\mbox{$\tilde{\bf E}_{3}$}}
\newcommand{\eiii}{\mbox{${\bf E}_{3}$}}
\newcommand{\eii}{\mbox{${\bf E}_{2}$}}
\newcommand{\ei}{\mbox{${\bf E}_{1}$}}
\newcommand{\eo}{\mbox{${\bf E}_{0}$}}
\newcommand{\ephi}{\mbox{$(E, \phi)$}}
\newcommand{\et}{\tilde{E}}
\newcommand{\fb}{\mbox{$\bf F$}}
\newcommand{\fp}{\mbox{     $\Box$}}
\newcommand{\fo}{\mbox{${\bf F}_{0}$}}
\newcommand{\gie}{\cite{g}}
\newcommand{\gliib}{\mbox{$\overline{\rm GL}$}}
\newcommand{\half}{{\scriptstyle{\frac{1}{2}}}}
\newcommand{\hx}{\raisebox{.5ex}{$\x$}\!}
\newcommand{\iione}{\mbox{${\rm II}_{1}$}}
\newcommand{\iitwo}{\mbox{${\rm II}_{2}$}}
\newcommand{\ip}{p_{1}}
\newcommand{\is}{\mbox{${\rm I}_{s}$}}
\newcommand{\iu}{\mbox{${\rm I}_{u}$}}
\newcommand{\ix}{\mbox{$X_{1}$}}
\newcommand{\japi}{\mbox{$J_{\alpha + 1}$}}
\newcommand{\ja}{\mbox{$J_{\alpha}$}}
\newcommand{\lat}{\mbox{$\tilde{\Lambda}$}}
\newcommand{\lii}{\Lambda^{2}}
\newcommand{\lli}{\Lambda L^{-1}}
\newcommand{\lj}{\lll_{j}}
\newcommand{\lll}{m}
\newcommand{\lot}{\mbox{$\tilde{\Lambda}_{0}$}}
\newcommand{\low}{\mbox{$\leftarrow$}}
\newcommand{\lrow}{\mbox{$\longrightarrow$}}
\newcommand{\lo}{\mbox{$\Lambda_{0}$}}
\newcommand{\lt}{\mbox{$\tilde{L}$}}
\newcommand{\mbt}{\tilde{M}_b}
\newcommand{\miotm}{\tilde{M}^-_{i-1}}
\newcommand{\mitl}{\tilde{M}_i}
\newcommand{\mitp}{\tilde{M}^+_i}
\newcommand{\mj}{\mmm_{j}}
\newcommand{\ml}{\mbox{$M^{\Lambda}$}}
\newcommand{\mmm}{n}
\newcommand{\mt}{\mbox{$\tilde{M}$}}
\newcommand{\nl}{\mbox{$N^{\Lambda}$}}
\newcommand{\nt}{\mbox{$\tilde{N}$}}
\newcommand{\pf}{{\em Proof}}
\newcommand{\pii}{p_{2}}
\newcommand{\qt}{\mbox{$\tilde{Q}$}}
\newcommand{\rank}{\mathop{{\rm rank}}\nolimits}
\newcommand{\row}{\mbox{$\rightarrow$}}
\newcommand{\s}{\sigma}
\newcommand{\siii}{\mbox{$S_{3}$}}
\newcommand{\siv}{\mbox{$S_{4}$}}
\newcommand{\slb}{\mbox{$\overline{\rm SL}$}}
\newcommand{\sliib}{\mbox{$\overline{\rm SL}_{2}$}}
\newcommand{\sliic}{\mbox{${\rm SL}_{2}{\rm(\C)}$}}
\newcommand{\slii}{\mbox{${\rm SL}_{2}$}}
\newcommand{\slx}{\mbox{SL($\x$)}}
\newcommand{\so}{\mbox{$S_{0}$}}
\newcommand{\sut}{\mbox{${\rm SU(2)}$}}
\newcommand{\td}{\mathop{{\rm td}}\nolimits}
\newcommand{\tensor}{\otimes}
\newcommand{\tr}{\mathop{{\rm tr}}\nolimits}
\newcommand{\ut}{\mbox{$\tilde{\cal U}$}}
\newcommand{\vo}{\mbox{${\cal W}^{\bla}_{0}$}}
\newcommand{\vol}{\mathop{{\rm vol}}\nolimits}
\newcommand{\vot}{\mbox{$\tilde{\cal W}^{\bla}_{0}$}}
\newcommand{\wl}{\mbox{${\cal W}^{\bla}$}}
\newcommand{\wot}{\mbox{$\tilde{\cal W}_{0}$}}
\newcommand{\wo}{\mbox{${\cal W}_{0}$}}
\newcommand{\wz}{\mbox{${\cal W}_{z}$}}
\newcommand{\x}{\chi}
\newcommand{\xii}{\mbox{$X_{2}$}}
\newcommand{\xn}{\mbox{$X_{n}$}}
\newcommand{\xo}{\mbox{$X_{0}$}}
\newcommand{\xt}{\mbox{$\tilde{X}$}}
\newcommand{\xz}{\mbox{$X_{z}$}}
\newcommand{\zk}{\mbox{$Z_{k}$}}
\newcommand{\zkl}{\mbox{$Z_{k}^{\Lambda}$}}
\newcommand{\zklt}{\mbox{$Z_{k}^{\tilde{\Lambda}}$}}
\newcommand{\zt}{\mbox{$\tilde{Z}$}}
\renewcommand{\det}{\mathop{{\rm det}}\nolimits}
\renewcommand{\bullet}{-}

\catcode`\@=12
\noindent
{\LARGE \bf Stable pairs, linear systems} \\
{\LARGE \bf and the Verlinde formula} \\
$\phantom{.}$ \\
{\bf
Michael Thaddeus }\\ \smallskip 
Mathematical Sciences Research Institute, 1000 Centennial Drive, Berkeley, Cal.
 94720 \\ 

\setcounter{section}{-1}
\bit{Introduction}

Let $X$ be a smooth projective complex curve of genus $g \geq 2$, let $\La \row
X$ be a line bundle of degree $d > 0$, and let \ephi\ be a pair consisting of a
vector bundle $E \row X$ such that $\lii E = \La$ and a section $\phi \in
H^0(E) - 0$.  This paper will study the moduli theory of such pairs.  However,
it is by no means a routine generalization of the well-known theory of stable
bundles.  Rather, it will discuss at least three remarkable features of the
moduli spaces of pairs:

1. Unlike bundles on curves, pairs admit many possible stability conditions.
In fact, stability of a pair depends on an auxiliary parameter $\s$ analogous
to the weights of a parabolic bundle.  This parameter was first detected by
Bradlow-Daskalopoulos \cite{bd} in the study of vortices on Riemann surfaces,
and indeed the spaces we shall construct can also be interpreted as moduli
spaces of rank 2 vortices.  As $\s$ varies, we will see that the moduli space
undergoes a sequence of flips in the sense of Mori theory, whose locations can
be specified quite precisely.

2. For some values of $\s$ the moduli space $M(\s, \La)$ is the blow-up of
$\Pj H^1(\La^{-1})$ along $X$, embedded as a complete linear system.  Thus we
can use $M(\s, \La)$ to study the projective embeddings of $X$.  In particular,
we obtain a very general formula \re{6m} for the dimension of the space of
hypersurfaces of degree $m+n$  in $\Pj H^1(\La^{-1})$ with a singularity at $X$
of order $n-1$.  This formula does not depend on the precise choice of $X$ and
$\La$, only on $g$ and $d$, which is rather surprising.

3.  For other values of $\s$, stability of the pair implies semistability of
the bundle, so $M(\s, \La)$ plays the role in rank 2 Brill-Noether theory of
the symmetric product in the usual case, and there is an Abel-Jacobi map from
$M(\s, \La)$ to the moduli space of semistable bundles. For large $d$ this is
generically a fibration, so we can use moduli spaces of pairs to study moduli
spaces of bundles.  In particular, we recover the known formulas for Poincar\'e
polynomials \cite{ab,hn} and Picard groups \cite{dn}; more spectacularly, we
prove, and generalize, the rank 2 Verlinde formula \re{6n} for both odd and
even degrees.

We will not fully discuss the many other fascinating aspects of the subject,
but we will briefly touch on one of them---the relation with Cremona
transformations and Bertram's work on secant varieties---in an appendix, \S8.
We hope to treat the relation with vortices and Yang-Mills-Higgs theory in a
later paper.

An outline of the other sections is as follows.  In \S1 we prove some basic
facts about pairs, in analogy with bundles.  Following Gieseker \cite{g1}, we
then use geometric invariant theory to construct the moduli space $M(\s, \La)$
of $\s$-semistable pairs, and a universal family over the stable points of
$M(\s, \La)$.  The choice of $\s$ corresponds to a choice of linearization for
our group action.  In \S2 we discuss the deformation theory of the moduli
problem.  In \S3 we show that the $M(\s, \La)$ are reduced, rational, and
smooth at the stable points.  We then show that as $\s$ varies, $M(\s, \La)$
undergoes a sequence of flips whose centres are symmetric products of $X$.  We
also define the rank 2 Abel-Jacobi map mentioned above.  In \S4 we calculate
the Poincar\'e polynomial of $M(\s, \La)$, and extract from it the
Harder-Narasimhan formula for the Poincar\'e polynomial of the moduli space of
rank 2 bundles of odd degree.

Thereafter we concentrate on studying the line bundles over $M(\s, \La)$, and
their spaces of sections.  In \S5 we compute the Picard group of $M(\s, \La)$,
and its ample cone. We explain how any section of a line bundle on $M(\s, \La)$
can be interpreted as a hypersurface in projective space, singular to some
order on an embedded $X$.  We also make the connection with the Verlinde vector
spaces.  Finally in \S\S6 and 7 we use the Riemann-Roch theorem to calculate
Euler characteristics of the line bundles in $M(\s, \La)$.  Combined with the
information from \S5, Kodaira vanishing, and some residue calculations which
were carried out by Don Zagier, this gives a formula for the dimensions of the
spaces of sections of line bundles on $M(\s, \La)$, under some mild hypotheses.
 We conclude by extracting the Verlinde formula from this.

For convenience we work over the complex numbers, but much of the paper should
be valid over any algebraically closed field: certainly \S\S1--3 and 5.
Kodaira vanishing is of course crucial in \S6, but the computation of the Euler
characteristics ought to make sense in general, if integral cohomology is
replaced with intersection theory.

A few notational habits should be mentioned: $X_i$ refers to the $i$th
symmetric product of $X$; $\pi$ denotes any obvious projection, such as
projection on one factor, or down from a blow-up; tensor products of vector
bundles are frequently indicated simply by juxtaposition; and likewise a
pullback such as $f^*L$ is often called just $L$.  Also, in \S3 and thereafter,
$M(\s, \La)$ is referred to simply as $M_i$, where $i$ depends on $\s$ in a
manner explained in \S3.  These conventions are not meant to be elliptical, but
to clean up what would otherwise be some very messy formulas.

We also make the following assumptions, which are explained in the text but are
repeated here for emphasis.  We always assume $g \geq 2$.  In the geometric
invariant theory construction of \S1, we assume $d$ is large, an assumption
which is justified by \re{3o} and the discussion following it.  From \S3 to the
end we assume $d \geq 3$.  However, this assumption is implicit in other
inequalities---so for example our main formula \re{6m} is valid as it stands.

{\em Acknowledgments.}  My principal debt of gratitude is of course to Don
Zagier, whose exquisite computations are indispensable to the paper.  The proof
of \re{4u}, and the entire \S7, are due to him.  I am also very grateful to
Simon Donaldson for his advice, encouragement and patience, and to Arnaud
Beauville, Aaron Bertram, Steven Bradlow, Jack Evans, Oscar Garcia-Prada, Rob
Lazarsfeld, David Reed, Miles Reid, and Eve Simms for helpful conversations.
Finally, I thank Krzysztof Gaw\c edzki and the Institut des Hautes Etudes
Scientifiques for their hospitality while much of the research for this paper
was carried out.

\bit{Constructing moduli spaces of $\s$-semistable pairs}
Our main objects of study, which we refer to simply as {\em pairs}, will be
pairs \ephi\ consisting of a rank 2 algebraic vector bundle $E$ over our curve
$X$, and a nonzero section $\phi \in H^{0}(E)$.  A careful study of such pairs
was made by Steven Bradlow \cite{brad}.  He defined a stability condition for
pairs and proved a Narasimhan-Seshadri-type theorem relating stable pairs to
vortices on a Riemann surface.  The vortex equations depend on a positive real
parameter $\tau$, and so the stability condition also depends on $\tau$.
Bradlow and Georgios Daskalopoulos went on \cite{bd} to give a gauge-theoretic
construction of the moduli space of $\tau$-stable pairs, under certain
conditions on $\tau$ and $\deg E$.  Oscar Garcia-Prada later showed \cite{gp}
that there always exists a projective moduli space, by realizing it as a
subvariety of a moduli space of stable bundles on $X \times \Pj^1$.  In this
section we will give a geometric invariant theory construction of the moduli
space of $\tau$-stable pairs for arbitrary $\tau$ and $\deg E$ (though for
convenience we assume $\rank E = 2$).  Aaron Bertram has informed me that he
has done something similar \cite{bert2}, and I apologize to him for any
overlap.

The Bradlow-Daskalopoulos stability condition is in general rather complicated,
but in the rank 2 case it simplifies to the following.  Let $\s$ be a positive
rational number.  It is related to $\tau$ by $\s =  \tau \vol X/4 \pi - \deg E
/2$, where $\vol X$ is the volume of $X$ with respect to the metric chosen in
\cite{bd}.

\begin{defn}
\label{3l}
The pair \ephi\ is $\s$-{\em semistable} if for all line bundles $L \subset E$,
$$\begin{array}{cl}
\deg L \leq \half \deg E - \s & \mbox{if $\phi \in H^{0}(L)$ and} \\
\deg L \leq \half \deg E + \s & \mbox{if $\phi \not\in H^{0}(L)$.}
\end{array}$$
It is  $\s$-{\em stable} if both inequalities are strict.
\end{defn}
The main result of this section is then the following.

\begin{thm}
\label{3b}
Let $\La \row X$ be a line bundle of degree $d$.  There is a projective moduli
space $M(\s, \La)$ of $\s$-semistable pairs \ephi\ such that $\lii E = \La$,
nonempty if and only if $\s \leq d/2$.
\end{thm}

Our construction will be modelled on that of Gieseker \cite{g1}.  We begin with
a few basic facts about $\s$-stable and semistable pairs, parallel to those for
bundles.  We write $\La$ for $\lii E$, and $d$ for $\deg E = \deg \La$.

\begin{lemma}
For $\s > 0$, there exists a $\s$-semistable pair of determinant $\La$ if and
only if $\s \leq d/2$.
\end{lemma}

\pf.  If $\s > d/2$, then $\s$-semistability implies $\deg L < 0$ if $\phi \in
H^0(L)$, which is absurd.  If $\s \leq d/2$, let $L \row X$ be a line bundle of
degree $[d/2-\s]$ having a nonzero section $\phi$.  Let $E$ be a nonsplit
extension
$$0 \lrow L \lrow E \lrow \lli \lrow 0.$$
Then the first inequality in the definition \re{3l} is obvious.  As for the
second, if $M \subset E$ and $\deg M > d/2 + \s$, then there is a nonzero map
$M \row \lli$.  Since $\deg \lli < d/2 + \s + 1$, this is an isomorphism, so
the extension is split, which is a contradiction.  \fp

\begin{lemma}
\label{3a}
Let \ephi\ be a pair.  There is at most one $\s$-destabilizing bundle $L
\subset E$ such that $\phi \in H^{0}(L)$, and at most one $\s$-destabilizing $M
\subset E$ such that $\phi \not\in H^{0}(M)$.  If both $L$ and $M$ exist, then
$E = L \oplus M$.
\end{lemma}

\pf.  The first statement is obvious, and the second follows from the
uniqueness of ordinary destabilizing bundles, since $\deg M \geq \half \deg E +
\s > \half \deg E$.  If both $L$ and $M$ exist, then the map $M \row E \row
\lli$ is nonzero since $\phi \in H^{0}(L)$ but $\not\in H^{0}(M)$.  But $\deg M
\geq d/2 + \s \geq \deg \lli$, so $M = \lli$ and $E$ is split. \fp

\begin{lemma}
Let $(E_{1}, \phi_{1})$ and  $(E_{2}, \phi_{2})$ be $\s$-stable pairs of degree
$d$, and let $\psi: E_{1} \row E_{2}$ be a map such that $\psi \phi_{1} =
\phi_{2}$.  Then $\psi$ is an isomorphism.
\end{lemma}

\pf.  The kernel of $\psi$ is a subsheaf of a locally free sheaf on a smooth
curve, so it is locally free.  If $\rank\ker \psi = 2$, then $\psi$ is
generically zero, so $\psi = 0$ and $\psi \phi_{1} \neq \phi_{2}$.  If
$\rank\ker \psi = 1$, then $\ker \psi$ is a line subbundle $L$ of $E_{1}$,
since $E_{1}/\ker \psi$ is contained in the torsion-free sheaf $E_{2}$.  Hence
$\psi$ descends to a map $\lli \row E_{2}$ (possibly with zeroes) such that
$\phi_{2} \in H^{0}(\lli)$.  Since $(E_{2}, \phi_{2})$ is $\s$-stable, $\deg
\lli < d/2 - \s$, so $\deg L > d/2 + \s$, contradicting the $\s$-stability of
$(E_{1}, \phi_{1})$.  Finally, if  $\rank\ker \psi = 0$, then $\ker \psi = 0$
and $\psi$ is injective.  Moreover, $\coker\psi$ is a coherent sheaf on a curve
with rank and degree 0, so $\coker\psi = 0$ and $\psi$ is an isomorphism. \fp

\begin{lemma}
\label{3k}
Let \ephi\ be a $\s$-stable pair. Then there are no endomorphisms of $E$
annihilating $\phi$ except 0, and no endomorphisms preserving $\phi$ except the
identity.
\end{lemma}

\pf. Subtracting from the identity interchanges the two statements, so they are
equivalent.  We prove the second.  Any endomorphism annihilating $\phi$
annihilates the subbundle $L$ generated by $\phi$, so descends to a map $E/L
\row E$.  But by $\s$-stability $E/L$ is a line bundle of degree $\geq d/2 +
\s$, so the image of this map, if it were nonzero, would generate a line bundle
of degree $\geq d/2 + \s$, which would be destabilizing.  \fp

\begin{lemma}
\label{3n}
Let $(\be,\bp), (\be',\bp') \row T \times X$ be two families over $T$
parametrizing the same pairs.  Then $(\be,\bp) \conga (\be',\bp')$.
\end{lemma}

\pf.  For any $t \in T$, the subspace of $H^0(X;\Hom(\be_t,\be'_t))$ consisting
of homomorphisms $\psi$ such that $\psi \bp_t = \lambda \bp'_t$ for some
$\lambda \in \C$ is one-dimensional by \re{3k}.  This determines an invertible
subsheaf of the direct image $(R^0\pi)\Hom(\be_t,\be'_t)$.  But this subsheaf
is trivialized by the section $\lambda = 1$, which produces the required
isomorphism.  \fp

The notion of a Harder-Narasimhan filtration for rank 2 pairs is quite a simple
one.  For \ephi\ stable, define $\Gr\ephi = \ephi$.  Otherwise, define
$\Gr\ephi$ to be a direct sum of line bundles, one of them containing the
section $\phi$, as follows.  If $L$ is the destabilizing bundle and $\phi \in
H^{0}(L)$, define $\Gr\ephi = (L \oplus \lli, \phi)$.   If $M$ is the
destabilizing bundle and $\phi \not\in H^{0}(M)$, project $\phi$ to a nonzero
section $\phi' \in H^{0}(\La M^{-1})$ and define $\Gr\ephi = (M \oplus \La
M^{-1}, \phi')$.  Note that if there are destabilizing bundles of both sorts,
then by \re{3a} $E = L \oplus \lli$ and the two definitions agree.

\begin{lemma}
\label{3j}
There exists a degeneration of \ephi\ to $\Gr\ephi$, but $\Gr\ephi$ degenerates
to no semistable bundle.
\end{lemma}

\pf.  The first statement is vacuous when \ephi\ is stable.  If it is unstable,
say with destabilizing bundle $M$, we can construct a pair $(\be, \bp) \row X
\times \C$ such that $(\be_{z}, \bp_{z}) \cong \ephi$ for $z \neq 0$, but
$(\be_{0}, \bp_{0}) \cong \Gr\ephi$, as follows.  Pull back \ephi\ to $X \times
\C$, and tensor by $\co(0)$ when $\phi \not\in H^0(M)$.  This gives a pair
$(\be', \bp') \row X \times \C$ such that $\bp'$ is annihilated by the natural
map $\be' \row \La M^{-1}|_{X \times \{ 0 \} }$.  Let $\be$ be the kernel of
this map; then $\bp'$ descends to $\bp \in H^0(\be)$, and it is straightforward
to check that $(\be, \bp)$ has the desired properties.

As for the second statement, suppose first that \ephi\ is stable.  If $C$ is a
curve, $p \in C$, and $(\be, \bp) \row X \times C$ is a flat family of pairs
such that $(\be_{z}, \bp_{z}) \cong \ephi$ for $z \neq p$, then $\bp_{p}$ has
the same zero-set $D$ as $\phi$, so $E$ and $\be_{p}$ are both extensions of $L
= \co(D)$ by $\La(-D)$; indeed, $\be$ is a family of such extensions.  The
extension class varies continuously, so the extension class of $\be_{p}$ is in
the same ray as that of $E$.  If it is nonzero, $\ephi \cong (\be_{p},
\bp_{p})$, and if it is zero, $(\be_{p}, \bp_{p})$ is destabilized by $\lli$.

Now suppose that \ephi\ is not stable, so that for some $L$, $\Gr\ephi = L
\oplus \lli$ and $\phi \in H^{0}(L)$.  Then as above $\be_{p}$ is an extension
of $L$ by $\lli$, but now by continuity the extension class must be zero, so
$\Gr\ephi = (\be_{p}, \bp_{p})$. \fp

\begin{lemma}
\label{3o}
If \ephi\ is $\s$-(semi)stable, then so is $(E(D), \phi(D))$ for any effective
divisor $D$.  Likewise, if $\phi$ vanishes on an effective divisor $D$ and
\ephi\ is $\s$-(semi)stable, then so is $(E(-D), \phi(-D))$.
\end{lemma}

\pf.  If $L \subset E$ is any line bundle, $\phi(D) \in H^{0}(L(D))$ if and
only if $\phi \in H^{0}(L)$, and $\deg L(D) = \deg L + \deg D$.  But $\half
\deg E(D) = \half \deg E + \deg D$ also, so both inequalities are preserved by
tensoring with $D$.  The second statement is proved similarly.  \fp

Hence if the moduli spaces $M(\s, \La)$ exist for large enough $d$, then the
moduli spaces for smaller $d$ will be contained inside them as the locus of
pairs \ephi\ such that $\phi$ vanishes on some effective $D$.  So to prove our
existence theorem \re{3b} it suffices to construct $M(\s, \La)$ for $d$ large
relative to $g$ and $\s$, and {\em we will assume for the remainder of\/ {\rm
\S\thesection\ }that $d$ is large in this sense}.  For such a large $d$, we
then have the following useful fact.

\begin{lemma}
\label{3e}
For fixed $g$ and $\s$ and large $d$, \ephi\ $\s$-semistable implies that
$H^{1}(E) = 0$ and $E$ is globally generated.
\end{lemma}

\pf.  Suppose that $H^{1}(E) \neq 0$.  Then $H^{0}(KE^{*}) \neq 0$, so there is
an injection $0 \row K^{-1}(D) \row E^{*}$ for some effective $D$.  Hence there
is an injection $0 \row  K^{-1}\La(D) \row E$.  Since $\deg K^{-1}\La(D) \geq 2
- 2g + d$, the $\s$-semistability condition implies that $2 - 2g + d \leq d/2 +
\s$, so that $d \leq 4g-4+2\s$.  So for $d$ larger than this, $H^{1}(E) = 0$.

Similarly, if $d > 4g-2+2\s$, then $H^{1}(E(-x)) = 0$ for all $x \in X$, so $E$
is globally generated.  \fp

Since we are assuming that $d$ is large, the above lemma implies that for
\ephi\ $\s$-stable, $\dim H^{0}(E) = \x(E) = d+2-2g$.  Call this number $\x$.
If we fix an isomorphism $s: \C^{\x} \row H^{0}(E)$, we obtain a map $\lii
\C^{\x} \stackrel{s}{\row} \lii H^{0}(E) \stackrel{\wedge}{\row} H^{0}(\La)$,
which is nonzero because $E$ is globally generated.  Thus to any bundle $E$
appearing in a $\s$-semistable pair, and any isomorphism $s$, we associate a
point $T(E,s) \in \Pj \Hom(\lii \C^{\x}, H^{0}(\La))$.  We will consider the
pair $(T(E,s), s^{-1}\phi) \in \Pj \Hom \times \Pj\C^{\x}$, where $\Pj \Hom$ is
short for $\Pj \Hom(\lii \C^{\x}, H^{0}(\La))$.  Roughly speaking, $M(\s, \La)$
will be a geometric invariant theory quotient of the set of such pairs.  The
quotient is necessary to remove the dependence on the choice of $s$.  Since two
such isomorphisms are related by an element of $\slx$, the group action will be
the obvious diagonal action of $\slx$ on $\Pj \Hom \times \Pj\C^{\x}$.  As
usual in geometric invariant theory, we must {\em linearize} the action by
choosing an ample line bundle and lifting the action of $\slx$ to its dual.  So
let the ample bundle be any power of $\co(\x + 2\s, 4\s)$, with the obvious
lifting.  (Of course $\x + 2\s$ and $4\s$ may not be integers, but by abuse of
notation we will refrain from clearing denominators, since the choice of power
does not matter.)  We can then define stable and semistable points in the sense
of geometric invariant theory with respect to this linearization.

\begin{propn}
\label{3f}
If \ephi\ is $\s$-(semi)stable, then $(T(E,s), s^{-1}\phi)$ is a (semi)stable
point with respect to the linearization above.
\end{propn}

\pf.  Suppose $T = (T(E,s), s^{-1}\phi)$ is not semistable.  Then by Mumford's
numerical criterion \cite{mf,new} there exists a nontrivial 1-parameter
subgroup $\lambda: \C^{\times} \row \slx$ such that for any $\tilde{T}$ in the
fibre of the dual of our ample bundle over $T$, $\lim_{t \row 0}\lambda(t)
\cdot \tilde{T} = 0$.  We interpret this limit concretely as follows.  Any
1-parameter subgroup of $\slx$ can be diagonalized, so there exists a basis
$e_{i}$ of $\C^{\x}$ such that $\lambda(t) \cdot e_{i} = t^{r_{i}}e_{i}$, where
$r_{i} \in \Z$ are not all zero and satisfy $\sum_{i} r_{i} = 0$ and $r_{i}
\leq r_{j}$ for $i \leq j$.  Then $\lim_{t \row 0}\lambda(t) \cdot \tilde{T} =
0$ means that any basis element $(e_i^* \wedge e_j^* \otimes v, e_k) \in
\Hom(\lii \C^{\x}, H^0(\La)) \oplus \C^{\x}$ which is acted on with weight
$\leq 0$ has coefficient zero in the basis expansion of $\tilde{T}$.  Because
of our choice of linearization, this means that $T(E,s)(e_{i},e_{j}) = 0$
whenever
\begin{equation}
\label{3c}
r_{i} + r_{j} \leq \frac{2\s}{\hx/2 + \s}\,\, r_{\ell},
\end{equation}
where $\ell = \max \{ i: \mbox{coefficient of $e_{i}$ in $s^{-1}\phi$ is $\neq
0$} \}$.  Let $L \subset E$ be the line bundle generated by $s(e_{1})$.  We
distinguish between two cases, according to whether $\phi \in H^{0}(L)$.

First case: $\phi \in H^{0}(L)$.  For $i \leq \hx/2 - \s + 1$, note that
$$ (\hx/2 - \s)\, r_{1} + (\hx/2 + \s)\, r_{i} \leq \sum_{i} r_{i} = 0,$$
since the left-hand side can be regarded as the integral over $[0,\x)$ of a
(two-step) step function whose value on $[j-1, j)$ is $\leq r_{j}$.  Hence for
$i \leq \hx/2 - \s + 1$,
$$ r_{1} + r_{i} \leq \frac{2\s}{\hx/2 + \s}\,\, r_{1} \leq \frac{2\s}{\hx/2 +
\s}\,\, r_{\ell},$$
so $T(E,s)(e_{1}, e_{i}) = s(e_{1}) \wedge s(e_{i}) = 0$.  Hence $s(e_{i})$ is
a section of the same line bundle as $s(e_{1})$, namely $L$.  So $\dim H^{0}(L)
> \hx/2 - \s$; since $d$ is large relative to $g$ and $\s$, this implies that
$\deg L > d/2 - \s$, so \ephi\ is not $\s$-semistable.

Second case: $\phi \not\in H^{0}(L)$.  For $i \leq \hx/2 + \s + 1$,
$$ (\hx/2 + \s)\, r_{1} + (\hx/2 - \s)\, r_{i} \leq 0,$$
for the same reason as above.  Hence
$$ r_{1} + r_{i} \leq \frac{2\s}{\hx/2 + \s}\,\, r_{i}.$$
We claim that $\ell > \hx/2 + \s + 1$.  If not, then for all $i \leq \ell$,
$$r_{1} + r_{i} \leq \frac{2\s}{\hx/2 + \s}\,\, r_{\ell},$$
so that $s(e_{i})$ would be in the same line bundle as $s(e_{1})$.  Since
$\phi$ is a linear combination of $e_{i}$ for $i \leq \ell$, we would conclude
$\phi \in H^{0}(L)$, a contradiction.  This proves the claim.

So for $i \leq \hx/2 + \s + 1$, actually
$$ r_{1} + r_{i} \leq \frac{2\s}{\hx/2 + \s}\,\, r_{\ell};$$
hence $s(e_{i}) \in H^{0}(L)$ as in the first case.  So $\dim H^{0}(L) > \hx/2
+ \s$, and again \ephi\ is not $\s$-semistable.

The proof for stability is similar: the numerical criterion now just says
$\lim_{t \row 0}\lambda(t) \cdot \tilde{T} \neq \infty$, so we replace the
$\leq$ in \re{3c} by $<$.  We just need to note that if $i < \hx/2 - \s + 1$,
then
$$ (\hx/2 - \s)\, r_{1} + (\hx/2 + \s)\, r_{i} < 0$$
strictly, because either the two step functions are different just to the left
of $\hx/2 - \s$, or the smaller one is identically $r_{1} < 0$.  \fp

\begin{propn}
\label{3g}
Let \ephi\ be a pair, let $s: \C^{\x} \row H^{0}(E)$ be a linear map, and let
$v \in \C^{\x}$ satisfy $s(v) = \phi$.  Write $T_{s}$ for the composition $\lii
\C^{\x} \stackrel{s}{\row} \lii H^{0}(E) \stackrel{\wedge}{\row} H^{0}(\La)$.
If $(T_{s}, v)$ is semistable, then $s$ is an isomorphism and \ephi\ is
$\s$-semistable.
\end{propn}

\pf.  First of all, note that if $s$ is not injective, then $(T_{s}, v)$ is
certainly not semistable.  Indeed, if $s(w) = 0$ for some $w$,  put $e_{1} =
w$, $e_{2} = v$, extend to a basis $\{ e_{i} \}$ of $\C^{\x}$, and then take
the 1-parameter subgroup defined by $r_{1} = -\x + 2$, $r_{2} = 0$, $r_{3} =
\cdots = r_{\x} = 1$.  Then $\ell = 2$, so
$$r_{i} + r_{j} \leq \frac{2\s}{\hx/2 + \s}\,\, r_{l}$$
means just $r_{i} + r_{j} \leq 0$.  Hence either $i=1$, or $j=1$, or $i=j=2$;
in any case, clearly $T_{s}(e_{i},e_{j}) = 0$.

Suppose then that $s$ is injective and \ephi\ is $\s$-unstable.  We will prove
$(T_{s}, v)$ is unstable.  Let $L \subset E$ be the destabilizing bundle.  We
distinguish two cases, depending on the sign of $d - \deg L - 2g + 2$.

First case: $d - \deg L > 2g-2$.  Then $H^{1}(\lli) = 0$, but $H^{1}(L) = 0$
also since $\deg L > d/2 - \s$ which is large relative to $g$.  Hence from the
long exact sequence of
\beq
\label{3d}
0 \lrow L \lrow E \lrow \lli \lrow 0
\eeq
we find that $H^{1}(E) = 0$, so $\dim H^{0}(E) = \x$ and $s$ is an isomorphism.
 Choose a basis $ e_{1}, \dots , e_{p}$ for $s^{-1}(H^{0}(L))$ and extend to a
basis $ e_{1}, \dots , e_{\x}$ for $\C^{\x}$.  Take the 1-parameter subgroup
defined by $r_{i} = p-\x$ for $i \leq p$, $p$ for $i > p$.  Then $r_{\ell} = p
- \x$ if $\phi \in H^{0}(L)$, $p$ if $\phi \not\in H^{0}(L)$.  Since $L$ is
destabilizing, $p > \hx/2 - \s$ if $\phi \in H^{0}(L)$, $p > \hx/2 + \s$ if
$\phi \not\in H^{0}(L)$.
Either way,
$$r_{i} + r_{j} \leq \frac{2\s}{\hx/2 + \s}\,\, r_{l}$$
implies $i,j \leq p$; if $\phi \in H^{0}(L)$, and say $i > p$, then
$$r_{i} + r_{j} - \frac{2\s}{\hx/2 + \s}\,\, r_{l} \geq p + (p-\x)(1 -
\frac{2\s}{\hx/2 + \s}) = p\,\,\frac{\x}{\hx/2+\s} -
\x\,\frac{\hx/2-\s}{\hx/2+\s}$$
$$> (\hx/2 - \s)\frac{\x}{\hx/2+\s} - \x\,\frac{\hx/2-\s}{\hx/2+\s} = 0,$$
whereas if $\phi \not\in H^{0}(L)$, and say $j > p$, then
$$r_{i} + r_{j} - \frac{2\s}{\hx/2 + \s}\,\, r_{l} \geq p - \x + p\,(1 -
\frac{2\s}{\hx/2 + \s}) = p\,\,\frac{\x}{\hx/2+\s} - \x > \x - \x = 0.$$
But if $i, j \leq p$, then $s(e_{i}), s(e_{j}) \in H^{0}(L)$, so $T_{s}(e_{i},
e_{j}) = 0$.  Hence $(T_{s}, v)$ is unstable.

Second case: $d - \deg L \leq 2g-2$.  Then $\dim H^{0}(\lli) \leq g$, so from
the long exact sequence of \re{3d} we deduce that the codimension of $H^{0}(L)$
in $H^{0}(E)$ is $\leq g$.  Hence the codimension of $s^{-1}(H^{0}(L))$ in
$\C^{\x}$ is $\leq g$.  Choose a basis $ e_{1}, \dots , e_{p}$ for
$s^{-1}(H^{0}(L))$ and extend to a basis $ e_{1}, \dots , e_{\x}$ for
$\C^{\x}$.  Take the 1-parameter subgroup defined by $r_{i} = p-\x$ for $i \leq
p$, $p$ for $i > p$.  Since $p \geq \x - g$ and $\x = d + 2 - 2g$ is large
relative to $\s$ and $g$, certainly $p > \hx/2 + \s$.  The remainder of the
proof proceeds as in the first case.

So far we have proved that if $(T_{s}, v)$ is semistable, then $s$ is injective
and \ephi\ is $\s$-semistable.  But then by \re{3e}, $\dim H^{0}(E) = \x$, so
$s$ is an isomorphism.  \fp

\begin{propn}
\label{3h}
Suppose $(E_{1}, \phi_{1})$ and $(E_{2}, \phi_{2})$ are $\s$-semistable, and
there exist $s_{1}, s_{2}$ such that $(T(E_{1}, s_{1}), s_{1}^{-1}\phi_{1}) =
(T(E_{2}, s_{2}), s_{2}^{-1}\phi_{2})$.  Then there is an isomorphism $(E_{1},
\phi_{1}) \cong (E_{2}, \phi_{2})$  under which $s_{1} \cong s_{2}$.
\end{propn}

\pf.  By \re{3e} each $E_{i}$ is globally generated, so the components
$s_{i}(e_{j}) \wedge s_{i}(e_{k})$ of $T(E_{i}, s_{i})$ give a map from $X$ to
the Grassmannian of $(\x - 2)$-planes in $\C^{\x}$ such that $E_{i}$ is the
pullback of the tautological rank 2 bundle, $\phi_{i}$ is the pullback of the
section defined by $s_{i}^{-1}(\phi_{i})$, and $s_{i}$ is the natural map from
$\C^{\x}$ to the space of sections of the tautological bundle.  So we can
recover $(E_{i}, \phi_{i})$ and $s_{i}$, up to isomorphism, from $(T(E_{i},
s_{i}), s_{i}^{-1}\phi_{i})$. \fp

\begin{propn}
\label{3i}
Let $C$ be a smooth affine curve and $p \in C$.  Let $(\be, \bp)$ be a locally
free family of pairs on $X \times C - \{ p \}$, and suppose \be\ is generated
by finitely many sections $s_{i}$.  Then after possibly rescaling $\bp$ by a
function on $C - \{ p \}$, $(\be, \bp)$ and the $s_{i}$ extend over $p$ so that
\be\ is still locally free, $\bp_{p} \neq 0$, and the $s_{i}$ generate
$\be_{p}$ at the generic point.
\end{propn}

The reason for proving the last fact is to ensure that $T(E,s)$ is nonzero at
$p$, so defines an element of $\Pj \Hom$.

\pf.  Choose an ample line bundle $L$ on $X \times C - \{ p \}$ such that
$\be^{*} \otimes L$ is globally generated.  Then \be\ embeds in a direct sum of
copies of $L$, and $\oplus_{j} L$ can be extended over $p$ as a sum of line
bundles in such a way that the $s_{i}$ extend too.  Consider the subsheaf of
the extended $\oplus_{j} L$ generated by the $s_{i}$.  This is a subsheaf of a
locally free sheaf, so it is torsion-free, and hence \cite{oss} has singular
set $S$ of codimension $\geq 2$.  Furthermore, it injects into its double dual,
whose singular set has codimension $\geq 3$ \cite{oss}, hence is empty.  Hence
the double dual is a locally free extension of \be\ over $p$, and is generated
by $s_{i}$ away from $S$.  As for $\bp$, it certainly extends with a possible
pole at $p$, so it is just necessary to multiply it by a function on $C$
vanishing to some order at $p$.  \fp

We can finally proceed to construct the geometric invariant theory quotient.
Consider the Grothendieck Quot scheme \cite{grot} parametrizing flat quotients
of $\co^{\x}_{X}$ with degree $d$, let $\Quot(\La) \subset \Quot$ be the
locally closed subset consisting of locally free quotients $E$ with $\lii E =
\La$, and let $U \subset \Quot(\La)$ be the open set where the quotient induces
an isomorphism $s: \C^{\x} \row H^{0}(E)$.  Then the pair $E,s$ specifies a
point in $U$.  By \re{3e}, if $(E_{p}, \phi)$ is $\s$-semistable for any
section $\phi$, then $p \in U$.

Now $U$ is acted upon by \slx\ in the obvious way, and the map $$T \times 1: U
\times \Pj \C^{\x} \row \Pj \Hom \times \Pj \C^{\x}$$ intertwines the group
actions on the two sets.  By \re{3f} and \re{3g}, the $\s$-semistable set
$V(\s) \subset U \times \Pj \C^{\x}$ is the inverse image of the semistable set
$V'(\s) \subset \Pj \Hom \times \Pj \C^{\x}$ with respect to the linearization
$\co(\x+2\s, 4\s)$.  In future, we restrict $T \times 1$ to a map $V(\s) \row
V'(\s)$.

Now Gieseker proves the following.

\begin{lemma}
Let $G$ be a reductive group and $M_{1}$ and $M_{2}$ be two $G$-spaces.
Suppose that $f: M_{1} \row M_{2}$ is a finite $G$-morphism and that a good
quotient $M_{2} / / G$ exists.  Then a good quotient $M_{1} / / G$ exists, and
the induced morphism $M_{1} / / G \row M_{2} / / G$ is finite. \fp
\end{lemma}

So to show that $V(\s)$ has a good quotient it suffices to prove:

\begin{lemma}
On $V(\s)$, $T \times 1$ is finite.
\end{lemma}

\pf.  By \re{3h}, $T \times 1$ is injective.  We use the valuative criterion to
check that $T \times 1$ is proper.  Let $C$ be a smooth curve, $p \in C$, and
let $\Psi: C - \{ p \} \row V(\s)$ be a map such that $(T \times 1) \Psi$
extends to a map $C \row V'(\s)$.  On $C - \{ p \}$, we then have a family
$(\be, \bp)$ of pairs such that \be\ is generated by the sections $s(e_{1}),
\dots , s(e_{\x})$.  By \re{3i}, on an open affine of $C$ containing $p$,
$(\be, \bp)$ extends over $p$ in such a way that $\bp_{p} \neq 0$ and the
$s(e_{i})$ generically generate $\be_{p}$.  Thus $T(\be_{p}, s)$ is defined,
and so by continuity $(T(\be_{p}, s), s^{-1}\bp_{p}) = ((T \times 1) \Psi)(p)$.
 Hence by \re{3g} $s: \C^{\x} \row H^{0}(\be_{p})$ is an isomorphism and
$(\be_{p}, \bp_{p})$ is $\s$-semistable.  So $(\be_{p}, s^{-1}\bp_{p}) \in
V(\s)$ and $\Psi$ extends to a map $C \row V(\s)$.  \fp

Hence $V(\s)$ has a good projective quotient.  By \re{3j}, the closure of the
orbit of $\ephi$ contains the orbit of $\Gr \ephi$, which is closed in the
$\s$-semistable set.  But the closure of any orbit in the $\x$-semistable set
contains only one closed orbit \cite[3.14 (iii)]{new}.  Hence if two pairs are
$\s$-semistable, then the closures of their orbits intersect if and only if
they have the same $\Gr$.  This completes the proof of our main theorem
\re{3b}. \fp

\begin{remk}
\label{3p}
If $D$ is any effective divisor, by \re{3o} there is an inclusion $\iota_D:
M(\s, \La) \hookrightarrow M(\s, \La(2D))$.  Indeed, if $(\be^{\La},
\bp^{\La})$ and $(\be^{\La(2D)}, \bp^{\La(2D)})$ are the corresponding
universal pairs, there is a sequence
$$0 \lrow \be^{\La} \stackrel{\iota_D}{\lrow} \iota_D^* \be^{\La(2D)} \lrow
\co_D (\iota_D^* \be^{\La(2D)}) \lrow 0$$
such that $\iota_D(\bp^{\La}) = \bp^{\La(2D)}$.
\end{remk}

One other pleasant fact should be mentioned: that the stable subsets of these
moduli spaces are fine.

\begin{propn}
\label{3m}
There exists a universal pair over the $\s$-stable set $M_s(\s,\La)$.
\end{propn}

\pf.  There is a universal bundle $\be \row \Quot(\La) \times X$ and a
surjective map $\co^{\x} \row \be$.  Hence there is a natural $\slx$-invariant
section $\bp \in H^0(\Quot(\La) \times \Pj\C^{\x} \times X; \be(1))$, and
$(\be(1), \bp)$ is a universal pair.  By \re{3k} the only stabilizers of
elements of the $\s$-stable subset of $V(\s)$ are the $\x$th roots of unity.
These act oppositely on $\be$ and on $\co(1)$, hence trivially on $\be(1)$, so
on the $\s$-stable set $\be(1)$ is invariant under stabilizers.  Hence by
Kempf's descent lemma \cite{dn} $\be(1)$ descends to a bundle on $M_s(\s, \La)
\times X$, and the section $\bp$, being invariant, also descends.  This pair
over $M_s(\s, \La) \times X$ then has the desired universal property.  \fp
\bit{Their tangent spaces}
We now turn to the deformation theory of our spaces.  By semicontinuity
$\s$-stability is an open condition, so the Zariski tangent spaces to our
moduli spaces at the $\s$-stable points will just be deformation spaces.  Hence
we may refer to $T_{(E, \phi)} M(\s, \La)$ simply as $T_{(E, \phi)}$.

\begin{propn}
\label{4c}
If $\ephi \in M(\s, \La)$ is $\s$-stable, then

{\rm (i) \,\,}(cf.\ {\rm \cite{bd}}) \,$T_{(E, \phi)}$ is canonically
isomorphic to $H^{1}$ of the complex
$$
C^{0}(\End_{0} E) \oplus \C \stackrel{p}{\lrow} C^{1}(\End_{0} E) \oplus
C^{0}(E) \stackrel{q}{\lrow} C^{1}(E),
$$
where $p(g,c) = (dg, (g+c)\phi)$ and $q(f, \psi) = f \phi - d\psi$;

{\rm (ii) \,\,}$H^0$ and $H^2$ of this complex vanish;

{\rm (iii) \,\,}there is a natural exact sequence
$$0 \lrow H^{0}(\End E) \stackrel{\phi}{\lrow} H^{0}(E) \lrow T_{(E, \phi)}
\lrow H^{1}(\End_{0} E) \stackrel{\phi}{\lrow} H^{1}(E) \lrow 0. $$
\end{propn}

\pf.  Let $R = \C[\varepsilon]/(\varepsilon^{2})$.  By a well-known result
\cite[II Ex.\ 2.8]{h} $T_{(E, \phi)}$ is the set of isomorphism classes of maps
$\Spec R \row M(\s, \La)$ such that $(\varepsilon) \mapsto \ephi$. Since
$\s$-stability is an open condition, $T_{(E, \phi)}$ is just the set of
isomorphism classes of families $(\be, \bp)$ of pairs on $X$ with base $\Spec
R$, such that $(\be, \bp)_{(\varepsilon)} = \ephi$ and $\lii \be$ is the
pullback of $\La$.  We will explain how to construct any such family.

The only open set in $\Spec R$ containing $(\varepsilon)$ is $\Spec R$ itself,
so any bundle $\be$ over $\Spec R \times X$ can be trivialized on $\Spec R
\times U_{\alpha}$ for some open cover $\{ U_{\alpha} \}$ of $X$.  Thus if
$\be_{(\varepsilon)} = E$, the transition functions give a \v{C}ech cochain of
the form $1 + \varepsilon f_{\alpha \beta}$ where $f \in C^{1}(\End E)$.  In
order for $\lii \be$ to be isomorphic to the pullback of $\La$, the transition
functions of $\lii \be$ must be conjugate to $1 \in C^{0}(\co)$.  But the
transition functions are $\det (1 + \varepsilon f_{\alpha \beta}) = 1 +
\varepsilon \tr f_{\alpha \beta}$, so we are asking that
$$(1 + \varepsilon g_{\alpha}) (1 + \varepsilon \tr f_{\alpha \beta}) (1 -
\varepsilon g_{\beta}) = 1$$
for some $g \in C^{0}(\co)$, that is, $\tr f = -d g$.  But if such a $g$
exists, then $\tilde{f} = f + dg/2$ is trace-free, and $1 + \varepsilon
\tilde{f}$ is obviously conjugate to $1 + \varepsilon f$, so determines the
same bundle \be.  Hence up to isomorphism we can obtain any \be\ even if we
consider only trace-free $f \in C^{1}(\End_{0} E)$.

Now if there is a section $\bp \in H^{0}(\be)$ such that $\bp_{(\varepsilon)} =
\phi$, then with respect to the local trivializations of \be\ described above,
$\bp = \phi + \varepsilon \psi_{\alpha}$ for some \v{C}ech cochain $\psi \in
C^{0}(E)$.  Of course, $\psi$ must be compatible with the transition functions;
this means that
$$(1 + \varepsilon f_{\alpha \beta}) (\phi + \varepsilon \psi_{\beta}) = (\phi
+ \varepsilon \psi_{\alpha}),$$ that is,  $f \phi = d \psi$.  Hence any pair
$(\be, \bp)$ having the desired properties can be obtained from some $(f, \psi)
\in C^{1}(\End_{0} E) \oplus C^{0}(E)$ satisfying $f \phi - d\psi = 0 \in
C^{1}(E)$.

We now need only check which $(f, \psi)$ give us isomorphic $(\be, \bp)$.  Of
course the two choices will be related by a change of trivialization on $\Spec
R \times U_{\alpha}$, but we may assume that the change of trivialization is of
the form $1 + \varepsilon g_{\alpha}$ on $U_{\alpha}$, since \ephi\ itself has
no automorphisms \re{3k}.  Furthermore, $g$ must belong to $C^{0}(\End_{0} E)
\oplus \C$ in order to keep $f$ trace-free, since the action of $g$ is given by
$$1 + \varepsilon f_{\alpha \beta} \mapsto (1 + \varepsilon g_{\alpha}) (1 +
\varepsilon f_{\alpha \beta}) (1 - \varepsilon g_{\beta}),$$
that is, $f \mapsto f + dg$, and $dg$ is trace-free if and only if $g \in
C^{0}(\End E)$ is the sum of a trace-free cocycle and a constant.  Similarly
the action of $g$ on $\psi$ is
$$\phi + \varepsilon \psi_{\alpha} \mapsto (1 + \varepsilon g_{\alpha}) (\phi +
\varepsilon \psi_{\alpha}), $$ that is, $\psi \mapsto \psi +  g \phi$.  Hence
two pairs $(f, \psi)$ and $(\tilde{f}, \tilde{\psi})$ determine isomorphic
pairs $(\be, \bp)$ if and only if they are in the same coset of the image of
the map $C^{0}(\End_{0} E) \oplus \C \row C^{1}(\End_{0} E) \oplus C^{0}(E)$
given by $g + c \mapsto (dg, (g+c)\phi)$. This completes the proof of (i).

As for (ii) and (iii), substituting $H^{0}(\End_{0} E) \oplus \C = H^{0}(\End
E)$ into the long exact sequence of the double complex with exact rows
$$
\begin{array}{ccccccccc}
0 & \lrow & 0 & \lrow & C^{0}(\End_{0} E) \oplus \C & \lrow &
C^{0}(\End_{0} E) \oplus \C & \lrow & 0 \\
 & & \Bdal{} & & \Bdal{} & & \Bdal{} & & \\
0 & \lrow & C^{0}(E) & \lrow & C^{1}(\End_{0} E) \oplus C^{0}(E) & \lrow &
C^{1}(\End_{0} E) & \lrow & 0 \\
 & & \Bdal{} & & \Bdal{} & & \Bdal{} & & \\
0 & \lrow & C^{1}(E) & \lrow & C^{1}(E) & \lrow & 0 & \lrow & 0
\end{array}
$$
gives
$$0 \lrow H^{0} \lrow H^{0}(\End E) \lrow  H^{0}(E) \lrow H^{1} \lrow
H^{1}(\End_{0} E) \lrow H^{1}(E) \lrow H^{2} \lrow 0,$$
where $H^{i}$ is the cohomology of the complex from (i). But the map
$H^{0}(\End E) \stackrel{\phi}{\lrow} H^{0}(E)$ is injective for \ephi\
$\s$-stable by \re{3k}, and the map $H^{1}(\End_{0} E) \stackrel{\phi}{\lrow}
H^{1}(E)$ is always surjective: indeed this is equivalent to the Serre dual map
$H^{0}(K E^* ) \stackrel{\phi}{\lrow} H^{0}(K \End_{0} E^{*})$ being injective,
which is obvious since the map $K E^{*} \stackrel{\phi}{\lrow} K \End_{0}
E^{*}$ is an injection of sheaves.  Hence $H^{0}$ and $H^{2}$ vanish, and we
get the exact sequence in (iii). \fp

As a corollary, we obtain the following.

\begin{cor}
\label{4w}
If $\ephi \in M(\s, \La)$ is $\s$-stable, then $\dim T_{(E, \phi)} = d+g-2$.
\end{cor}

\pf.  By \re{4c}(iii)
$$\dim T_{(E, \phi)} = \x(E) - \x(\End_{0} E) - 1 = (d+2-2g) - (3-3g) -1 =
d+g-2. \fp $$

We will see in the next section that $\dim M(\s, \La) = d+g-2$; hence $M(\s,
\La)$ will be smooth at the stable points.

\bit{How they vary with $\s$}

For obvious numerical reasons the $\s$-semistability condition remains the
same, and implies $\s$-stability, for any $\s \in (\max (0, d/2 - i - 1), d/2 -
i)$, where $i$ is an integer between 0 and $(d-1)/2$.  Hence for $\s$ in that
interval we get a fixed smooth projective moduli space $M(\s, \La)$, which we
will henceforth denote $M_{i}(\La)$ or just $M_{i}$.  The remainder of this
paper will concentrate on these smooth moduli spaces $M_i$, ignoring the
special values of $\s$ for which there exist $\s$-semistable pairs which are
not $\s$-stable.

In the extreme case $i=0$, it is then easy to construct the moduli space:

\beq
\label{4n}
M_{0}(\La) = \Pj H^{1}(\La^{-1}).
\eeq

\pf.  The first inequality in the $\s$-stability condition \re{3l} says that
$\phi \in H^{0}(L)$ implies $\deg L \leq 0$.  Hence $L = \co$, $E$ is an
extension of $\co$ by $\La$, and $\phi \in H^{0}(\co)$ is a constant section.
The second inequality says that $E$ has no subbundles of degree $\geq d$: this
is equivalent to not being split, since $M \row E \row \La$ nonzero and $\deg M
\geq d = \deg \La$ implies $M = \La$.  Hence $M_{0}(\La)$ is simply the moduli
space of nonsplit extensions of $\co$ by $\La$, which is of course just $\Pj
H^{1}(\La^{-1})$. \fp

We will not attempt such a direct construction of $M_{i}(\La)$ for $i > 0$.
Rather, we will carefully study the relationship between $M_{i-1}$ and $M_{i}$.
 Of course, this will only be of interest if there exists an $M_i$ for $i>0$,
so {\em we will assume for the remainder of the paper that $[(d-1)/2] > 0$,
that is, $d \geq 3$}.  Anyhow, the first step is to construct families
parametrizing those pairs which appear in $M_{i}$ but not $M_{i-1}$, or
$M_{i-1}$ but not $M_{i}$.  To do this, we first define two vector bundles over
the $i$th symmetric product $X_i$.

Let $\pi: X_i \times X \row X_i$ be the projection and let $\Delta \subset X_i
\times X$ be the universal divisor. Then define $W^-_i = (R^0 \pi)
\co_{\Delta}\La(-\Delta)$ and $W^+_i =  (R^{1}\pi)\La^{-1} (2\Delta)$.  These
are locally free sheaves of rank $i$ and $d+g-1-2i$, respectively.

\begin{propn}
\label{4g}
For $i \leq (d-1)/2$, there is a family over $\Pj W^+_i$ parametrizing exactly
those pairs which are represented in $M_i$ but not $M_{i-1}$.
\end{propn}

\pf. As we pass from $i$ to $i-1$, the first inequality in the stability
condition \re{3l} gets stronger and the second gets weaker.  So we look for
pairs which almost violate the first inequality.  That is, $E$ must be an
extension
$$0 \lrow \co(D) \lrow E \lrow \La(-D) \lrow 0,$$
where $\deg D = i$, and $\phi$ is the section of $\co(D)$ vanishing on $D$.
Conversely, any such pair is stable unless it splits $E = \co(D) \oplus
\La(-D)$.  Indeed, if $L \subset E$ and $\phi \not\in H^0(L)$, then the map $L
\row \La(-D)$ is nonzero, so $\deg L \leq \deg \La(-D) = d-i$, with equality
only if $L = \La(-D)$.

But $\Pj W^+_i$ is the base of a family parametrizing all such nonsplit pairs:
indeed \be\ is the tautological extension
$$0 \lrow \co(\Delta) \lrow \be \lrow \La(-\Delta)(-1) \lrow 0,$$
and $\bp$ is the section of $\co(\Delta)$ vanishing on $\Delta$.  \fp

\begin{propn}
\label{4h}
For $i \leq (d-1)/2$, there is a family over $\Pj W^-_i$ parametrizing exactly
those pairs which are represented in $M_{i-i}$ but not $M_{i}$.
\end{propn}

\pf.  This time the first inequality in \re{3l} gets weaker and the second gets
stronger.  So we look for pairs which almost violate the second inequality.
That is, $E$ is an extension
$$
0 \lrow M \lrow E \lrow \La M^{-1} \lrow 0
$$
where $\deg M = d-i$, and $\phi \not\in H^0(M)$.  Hence projecting $\phi$ in
the exact sequence, we get a nonzero $\gamma \in H^0(\La M^{-1})$ vanishing on
a divisor $D$ of degree $i$ such that $\La M^{-1} = \co(D)$.  Then at $D$,
$\phi$ lifts to $M = \La(-D)$, so we get an element $p\ephi \in
H^0(\co_D\La(-D))$, defined up to a scalar as usual.

On the other hand, we can recover \ephi\ from $D$ and $p$.  Indeed, choose a \v
Cech cochain $\psi \in C^0 (\La(-D))$ such that $\psi|_D = p$.  Then $d\psi|_D
= dp = 0$, so  $d\psi$ vanishes on $D$ and descends to a closed cochain $f =
d\psi/\gamma \in C^1((\La(-2D))$.  This determines an extension
$$ 0 \lrow \La(-D) \lrow E' \lrow \co(D) \lrow 0.$$
The compatibility condition for $\gamma + \psi$ to define a section $\phi' \in
H^0(E')$ is $\gamma f = d \psi$, which is automatic.  Thus we get a new pair
$(E', \phi')$ satisfying $p(E', \phi') = p$.

Up to isomorphism, $(E', \phi')$ is independent of the choice of $\psi$, since
adding $\xi \in C^0(\La(-2D))$ to $\psi$ is simply equivalent to acting by
$\left( \begin{array}{cc} 1 & \xi_{\alpha} \\ 0 & 1 \end{array} \right)$ on the
local splittings of $E'$ with which the extension is defined.  In particular,
we can choose local splittings of the old $E$ and let $\psi$ be the projection
of the old $\phi$ on $M = \La(-D)$ with respect to these splittings.  Then the
construction of the previous paragraph recovers \ephi, so $(E', \phi') \conga
\ephi$.

The construction above can be generalized to produce a family $(\be, \bp) \row
\Pj W^-_i \times X$, as follows.  Let $p: \Pj W^-_i \row X_i$ be the
projection, and choose a cochain $\Psi \in C^0(\La(-\Delta)(1))$ such that
$\Psi|_{p^{-1}\Delta}$ is the tautological section.  Then $d\Psi$ vanishes on
$p^{-1}\Delta$, so descends to $C^1(\La(-2\Delta)(1))$.  This determines an
extension
$$0 \lrow \La(-\Delta)(1) \lrow \be \lrow \co(\Delta) \lrow 0,$$
and if $\gamma \in H^0(\co(\Delta))$ is the section vanishing on $\Delta$, then
$\gamma + \Psi$ defines the desired section $\bp \in H^0(\be)$.  \fp

By the universal properties of $M_{i-1}$ and $M_i$, we thus get injections $\Pj
W^+_i \hookrightarrow M_{i}$ and $\Pj W^-_i \hookrightarrow M_{i-1}$.
As an example, consider the case $i = 1$.  By \re{4n}, $M_0 = \Pj
H^1(\La^{-1})$.  Moreover, $W^-_1$ is a line bundle and hence $\Pj W^-_1 = X_1
= X$.  Hence the inclusion of \re{4h} is a map $X \hookrightarrow \Pj
H^1(\La^{-1})$; it can be identified explicitly as follows.

\begin{propn}
The inclusion $X \hookrightarrow \Pj H^1(\La^{-1})$ is given by the complete
linear system $|K_X\La|$.
\end{propn}

\pf.  There is an alternative way to see what pairs are represented in $M_0$
but not $M_1$.  Any pair $\ephi \in M_0$ is an extension
\beq
\label{4o}
 0 \lrow \co \lrow E \lrow \La \lrow 0,
\eeq
say with extension class $s \in H^1(\La^{-1})$, and with $\phi \in H^0(\co)$.
Such a pair is the image of $x \in X$ under the injection of \re{4h} if there
is an inclusion $0 \row \La(-x) \row E$ such that the composition $\gamma_x:
\La(-x) \row E \row \La$ vanishes at $x$.  Hence we ask for what extension
classes $s \in H^1(\La^{-1})$ the map $\gamma_x: \La(-x) \row \La$ lifts to
$E$.

Twisting \re{4o} by $\La^{-1}(x)$ and taking the long exact sequence yields
$$H^0(E \otimes \La^{-1}(x)) \lrow H^0(\co(x)) \stackrel{s}{\lrow}
H^1(\La^{-1}(x)),$$
where the second map is the cup product with $s$.  Hence $\gamma_x \in
H^0(\co(x))$ lifts to $H^0(E \otimes \La^{-1}(x))$ as desired if and only if
$\gamma_x s = 0$.  That is, $s$ must be in the kernel of the map $\gamma_x:
H^1(\La^{-1}) \row H^1(\La^{-1}(x))$, or Serre dually, $\gamma_x: H^0(K_X\La)^*
\row H^0(K_X\La(-x))^*$.  Since $\gamma_x$ is dual to the injection
$H^0(K_X\La) \row H^0(K_X\La(-x))$, it is surjective, so
$$\dim \ker \gamma_x = \dim H^0(K_X\La(-x)) - \dim H^0(K_X\La).$$
But since $\deg K_X\La(-x) > 2g-2$, this is 1.  Hence for each $x \in X$, there
is a unique $s \in \Pj H^1(\La^{-1})$ such that $\gamma_x s = 0$.

What is this $s$?  Regarded as a linear functional on $H^0(K_X\La)$, $s \in
\ker \gamma_x$ if it annihilates all sections vanishing at $x$.  Certainly
evaluation at $x$ does this, so this is the $s$ generating $\ker \gamma_x$.
But it is also the image of $x$ in the map $X \hookrightarrow \Pj
H^1(\La^{-1})$ given by $|K_X\La|$.  Hence the two maps are identical.  \fp

\begin{propn}
\label{4t}
The $M_i$ are all smooth rational integral projective varieties of dimension $d
+g-2$, and for $i>0$, there is a birational map $M_i \leftrightarrow M_1$,
which is an isomorphism except on sets of codimension $\geq 2$.
\end{propn}

\pf.  By \re{4n} and Riemann-Roch, the first statement is certainly true of
$M_0$.  For $i>0$, suppose by induction on $i$ that it is true of $M_{i-1}$.
By \re{4g} and \re{4h} there is an isomorphism $M_{i-1} - \Pj W^-_i
\leftrightarrow M_i - \Pj W^+_i$.  But $\dim \Pj W^-_i = 2i-1 < d-1 < d+g-2$,
and $\dim \Pj W^+_i = d+g-2-i < d+g-2$, so $\dim M_i = \dim M_{i-1} = d+g-2$
and $M_i$ is birational to $M_{i-1}$, hence to $M_0$.  Moreover by \re{4w}, the
Zariski tangent space to $M_i$ has constant dimension $d+g-2$, so $M_i$ is a
smooth reduced variety.  The second statement is also proved by induction: just
note that for $i > 1$, $\codim \Pj W^-_i/M_{i-1} = d+g-2i-1 \geq 2$ and $\codim
\Pj W^+_i/M_i = i \geq 2$.  \fp

\begin{propn}
\label{4e}
Let $\ephi \in \Pj W^+_i$, let $D$ be the zero-set of $\phi$, and let $\gamma$
be the map $$E \otimes \La^{-1}(D) \row \La(-D) \otimes \La^{-1}(D) = \co.$$
Then
$T_{(E, \phi)}\Pj W^+_i$ is canonically isomorphic to $H^{1}$ of the complex
$$
C^{0}(E \otimes \La^{-1}(D)) \oplus \C \stackrel{p}{\lrow} C^{1}(E \otimes
\La^{-1}(D)) \oplus C^{0}(\co(D)) \stackrel{q}{\lrow} C^{1}(\co(D)),
$$
where $p(g,c) = (dg, (\gamma g + c) \phi)$ and $q(f, \psi) = \gamma f \phi -
d\psi$.  Moreover,
$H^0$ and $H^2$ of this complex vanish.
\end{propn}

\pf. The proof is modelled on that of \re{4c}.  We regard $\Pj W^+_i$ as a
moduli space of triples $(L, E ,\phi)$, where $L$ is a line bundle of degree
$i$, $E$ is an extension of $L$ by $\La L^{-1}$, and $\phi \in H^0(L)$, and
consider the deformation theory of this moduli problem.

Let $R = \C[\varepsilon ]/(\varepsilon ^2)$ as before. Then $T_{(L, E ,\phi)}
\Pj W^+_i$ is the set of isomorphism classes of families $({\bf L}, \be, \bp)$
of triples on $X$ with base $\Spec R$, such that $({\bf L}, \be,
\bp)_{(\varepsilon )} = (L, E ,\phi)$.  We will explain how to construct any
such family.

Any bundle over $\Spec R \times X$ can be trivialized on $\Spec R \times
U_{\alpha}$ for some open cover $\{ U_{\alpha} \}$ of $X$.  Thus if ${\bf
L}_{(\varepsilon )} = \co(D)$ and $\be_{(\varepsilon )} = E$, then the
transition functions for \be\ give a \v{C}ech cochain of the form $1 +
\varepsilon  f_{\alpha \beta}$ where $f \in C^{1}(\End E)$.  Since \be\ is to
be a family of extensions of $L$ by $\La L^{-1}$, it must have $\lii \be =
\La$, so as explained in the proof of \re{4c} we may take $f \in C^1(\End_0
E)$.  Furthermore, the transition functions must preserve $\bf L$, so if $f'$
is the projection of $f$ to $C^1(\La(-2D))$ in the natural exact sequence
$$0 \lrow E \otimes \La^{-1}(D) \lrow \End_0 E \lrow \La(-2D) \lrow 0,$$
then $1 + \varepsilon  f'_{\alpha \beta}$ must be conjugate to 1.  Hence
$$(1 - \varepsilon  g_{\alpha}) (1 + \varepsilon  f'_{\alpha \beta}) (1 -
\varepsilon  g_{\beta}) = 1$$
for some $g \in C^0(\La(-2D))$, that is, $f' = dg$.  But if such a $g$ exists,
then for any lifting $\tilde{g}$ of $g$ to $C^0(\End_0 E)$, $\tilde{f} = f -
d\tilde{g}$ projects to $0 \in C^1(\La(-2D))$, and $1 + \varepsilon  \tilde{f}$
is obviously conjugate to $1 + \varepsilon  f$, so determines the same bundle
\be.  Hence up to isomorphism we can obtain any \be\ that is an extension of
some $\bf L$ by $\La {\bf L}^{-1}$ even if we consider only those $f$ in the
kernel of $C^1(\End_0 E) \row C^1(\La(-2D))$, that is, in $C^1(E \otimes
\La^{-1}(D))$.  The transition functions for $\bf L$ are then just $1 +
\varepsilon  \gamma f_{\alpha \beta}$.

Now if there is a section $\bp \in H^0({\bf L})$ such that $\bp_{(\varepsilon
)} = \phi$, then with respect to the local trivializations of \be, $\bp = \phi
+ \varepsilon  \psi_{\alpha}$ for some \v{C}ech cochain $\psi \in C^0(\co(D))$.
Of course, $\psi$ must be compatible with the transition functions; this means
that
$$(1 + \varepsilon  \gamma f_{\alpha \beta}) (\phi + \varepsilon  \psi_{\beta})
= (\phi + \varepsilon  \psi_{\alpha}),$$ that is,  $\gamma f \phi = d \psi$.
Hence any triple $({\bf L}, \be, \bp)$ having the desired properties can be
obtained from some $(f, \psi) \in C^{1}(E \otimes \La^{-1}(D)) \oplus
C^{0}(\co(D))$ satisfying $\gamma f \phi - d\psi = 0 \in C^{1}(\co(D))$.

We now need only check which $(f, \psi)$ give us isomorphic $({\bf L},\be,
\bp)$.  This part of the argument follows that of \re{4c} exactly, except that
$g$ ends up being in $C^1(E \otimes \La^{-1}(D)) \oplus \C$, and acts on $\psi$
by $\psi \mapsto \psi + \gamma g \phi$.  This completes the proof of the first
statement.

As for the second, taking the long exact sequence of the double complex
$$
\begin{array}{ccccccccc}
0 & \lrow & 0 & \lrow & C^{0}(E \otimes \La^{-1}(D)) \oplus \C & \lrow &
C^{0}(E \otimes \La^{-1}(D)) \oplus \C & \lrow & 0 \\[.2em]
 & & \Bdal{} & & \Bdal{} & & \Bdal{} & & \\[.2em]
0 & \lrow & C^{0}(\co(D)) & \lrow & \def\arraystretch{.6}
\begin{array}{c} C^{1}(E \otimes \La^{-1}(D)) \\ \oplus C^{0}(\co(D))
\end{array} & \lrow &  C^{1}(E \otimes \La^{-1}(D)) & \lrow & 0 \\[.6em]
 & & \Bdal{} & & \Bdal{} & & \Bdal{} & & \\[.2em]
0 & \lrow & C^{1}(\co(D)) & \lrow & C^{1}(\co(D)) & \lrow & 0 & \lrow & 0
\end{array}
$$
gives
\beqas
\lefteqn{0 \lrow H^0 \lrow H^0(E \otimes \La^{-1}(D)) \oplus \C \lrow
H^0(\co(D)) \lrow H^1 } \\ & & \lrow H^1(E \otimes \La^{-1}(D)) \lrow
H^1(\co(D)) \lrow H^2 \lrow 0,
\eeqas
where $H^{i}$ is the cohomology of the complex in the statement.  Now
$H^0(\La^{-1}(2D)) = 0$ since $\deg \La^{-1}(2D) < 0$, and $E$ is a  nonsplit
extension of $\co(D)$ by $\La(-D)$, so
$$H^0(E \otimes \La^{-1}(D)) = H^0(\Hom(\La(-D), E)) = 0.$$
But the map $\C \row H^0(\co(D))$ is injective: indeed, it is multiplication by
$\phi$.  Hence $H^0 = 0$.  Likewise, the map $H^1(E \otimes \La^{-1}(D)) \row
H^1(\co(D))$ is surjective: indeed this is equivalent to the Serre dual map
$H^0(K(-D)) \row H^0(E^* \otimes K\La(-D))$ being injective, which is obvious
since the map $K(-D) \row K \row E^* \otimes K\La(-D)$ is an injection of
sheaves.  Hence $H^2 = 0$.  \fp

The following proposition is proved similarly.

\begin{propn}
\label{4v}
Let $\ephi \in \Pj W^-_i$, and let $D = p\ephi$.  Then
$T_{(E, \phi)}\Pj W^-_i$ is canonically isomorphic to $H^{1}$ of the complex
$$
C^0 (E(-D)) \oplus \C \lrow C^1 (E(-D)) \oplus C^0 (E) \lrow C^1 (E).
$$
Moreover, $H^0$ and $H^2$ of this complex vanish. \fp
\end{propn}

\begin{propn}
\label{4f}
The injection $\Pj W^+_i \hookrightarrow M_i$ induces an exact sequence on $\Pj
W^+_i$
$$0 \lrow T\Pj W^+_i \lrow TM_i |_{\Ps W^+_i} \lrow W^-_i(-1) \lrow 0.$$
\end{propn}

\pf.  The complex
$$
C^0 (\La(-2\Delta)) \lrow C^1 (\La(-2\Delta)) \oplus C^0 (\La(-\Delta)) \lrow
C^1 (\La(-\Delta))
$$
with the obvious maps has $R^0 \pi = 0$, $R^1 \pi = W^-_i$ from the long exact
sequence of the double complex
$$
\begin{array}{ccccccccc}
0 & \lrow & C^0 (\La(-2\Delta)) & \lrow & C^0(\La(-2\Delta)) & \lrow & 0 &
\lrow & 0 \\[.2em]
 & & \Bdar{(1,0)} & & \Bdar{} & & \Bdal{} & & \\[.2em]
0 & \lrow & \def\arraystretch{.6}
\begin{array}{c} C^0(\La(-2\Delta)) \\ \oplus C^1(\La(-2\Delta)) \end{array} &
\lrow & \def\arraystretch{.6}
\begin{array}{c} C^0(\La(-\Delta)) \\ \oplus C^1(\La(-2\Delta)) \end{array} &
\lrow &  C^0(\co_{\Delta} \La(-\Delta)) & \lrow & 0 \\[.6em]
 & & \Bdar{(0,1)} & & \Bdar{} & & \Bdal{} & & \\[.2em]
0 & \lrow & C^1(\La(-2\Delta)) & \lrow & C^1(\La(-\Delta)) & \lrow &
C^1(\co_{\Delta} \La(-\Delta))  & \lrow & 0.
\end{array}
$$
Hence the result follows from the long exact sequence of the double complex
$$
\begin{array}{ccccccccc}
0 & \lrow & 
 C^0 (\be \La^{-1}(\Delta)) \oplus \C
& \lrow & 
 C^0(\End_0 \be) \oplus \C
& \lrow & C^0 (\La(-2\Delta))(-1) & \lrow & 0 \\[.2em]
 & & \Bdar{(1,0)} & & \Bdar{p} & & \Bdal{} & & \\[.2em]
0 & \lrow &  \def\arraystretch{.6}
\begin{array}{c} C^1(\be \La^{-1}(\Delta)) \\ \oplus C^0(\co(\Delta))
\end{array}
& \lrow &  \def\arraystretch{.6}
\begin{array}{c} C^1(\End_0 \be) \\ \oplus C^0(\be) \end{array}
& \lrow &  \def\arraystretch{.6}
\begin{array}{c} C^1(\La(-2\Delta))(-1) \\ \oplus C^0(\La(-\Delta))(-1)
\end{array}
& \lrow & 0 \\[.6em] \def\arraystretch{1}
 & & \Bdar{(0,1)} & & \Bdar{q} & & \Bdal{} & & \\[.2em]
0 & \lrow & C^1(\co(\Delta)) & \lrow & C^1(\be) & \lrow & C^1(\La(-\Delta))(-1)
& \lrow & 0,
\end{array}
$$
together with \re{4c} and \re{4e}.  \fp

\begin{cor}
The map $\Pj W^+_i \hookrightarrow M_i$ is an embedding.
\end{cor}

\pf.  By \re{4e}, it is an injection, and by \re{4f}, so is its derivative.
\fp

The following proposition and corollary are proved similarly, using \re{4c} and
\re{4v}.

\begin{propn}
\label{4q}
The injection $\Pj W^-_i \hookrightarrow M_{i-1}$ induces an exact sequence on
$\Pj W^-_i$
$$0 \lrow T\Pj W^-_i \lrow TM_{i-1} |_{\Ps W^-_i} \lrow W^+_i(-1) \lrow 0.  \fp
$$
\end{propn}

\begin{cor}
The map $\Pj W^-_i \hookrightarrow M_{i-1}$ is an embedding.  \fp
\end{cor}

By \re{4g} and \re{4h} every pair in $M_i - \Pj W^+_i$ is also in $M_{i-1} -
\Pj W^-_i$, and vice-versa.  Hence there is a natural isomorphism $M_i - \Pj
W^+_i \row M_{i-1} - \Pj W^-_i$.  Our next task is to extend this to a proper
map.  Let $\mitp$ be the blow-up of $M_i$ at $\Pj W^+_i$.  Then by \re{4f}  the
exceptional divisor is $E^+_i = \Pj W^-_i \oplus \Pj W^+_i$, and
$\co_{E^+_i}(E^+_i) = \co(-1,-1)$.

\begin{propn}
\label{4l}
There is a map $\mitp \row M_{i-1}$ such that the following diagram commutes:
$$\begin{array}{ccccc}M_i - \Pj W^+_i & \lrow & \mitp & \longleftarrow & E^+_i
\\
                  \Big\updownarrow & & \Big\downarrow & & \Big\downarrow \\
                  M_{i-1} - \Pj W^-_i & \lrow & M_{i-1} & \longleftarrow & \Pj
W^-_i.\end{array}$$

\end{propn}

\pf.  Let $(\be, \bp) \row \mitp \times X$ be the pullback of the universal
family.  We will construct a new family $(\be', \bp')$ of pairs all of which
are in $M_{i-1}$.

By uniqueness of families \re{3n}, $(\be, \bp)|_{E^+_i \times X}$ is the
pullback of the family over $\Pj W^+_i$ constructed in \re{4g}.  Thus there is
a surjective sheaf map $\be \row  \co_{E^+_i \times X}\La(-\Delta)(0,-1)$
annihilating $\bp$.  Define $\be'$ to be the kernel of this map, so that
\beq
\label{4p}
0 \lrow \be' \lrow \be \lrow \co_{E^+_i \times X}\La(-\Delta)(0,-1) \lrow 0.
\eeq
Then $\be$ is locally free, and $\bp$ descends to $\bp' \in H^0(\be')$.  For $z
\in M_i - \Pj W^+_i$, clearly $(\be', \bp')_z = (\be, \bp)_z$.  So to prove the
proposition it suffices to show that $(\be', \bp')_{E^+_i}$ is the pullback of
the family over $\Pj W^-_i$ constructed in \re{4h}. The first  promising thing
to note is that there certainly a surjection $\be' \row \co_{E^+_i \times
X}(\Delta) \row 0$, and $\lii \be' = \lii \be(-E^+_i \times X)$, so we get an
extension
$$0 \lrow \La(-\Delta)(1,0) \lrow \be'_{E^+_i \times X} \lrow \co(\Delta) \lrow
0,$$
just as in the family of \re{4h}.

Now fix $s \in E^+_i$ over $\ephi \in M_i$, and let $D$ be the zero-set of
$\phi$.  Let $R = \C[\varepsilon]/(\varepsilon^2)$ as before, and choose a map
$\Spec R \row \mitp$ representing an element of $T_s \mitp - T_s E^+_i$.  Then
\re{4p} restricts to an exact sequence
$$0 \lrow \co_{\Spec R \times X}(\be') \lrow \co_{\Spec R \times X}(\be) \lrow
\co_{(\varepsilon) \times X}\La(-D) \lrow 0.$$
On some open cover $\{ U_{\alpha} \}$ of $X$, $E$ splits as
\beq
\label{4i}
E|_{U_{\alpha}} = \co(D)|_{U_{\alpha}} \oplus \La(-D)|_{U_{\alpha}},
\eeq
and this splitting can be extended to a splitting of $\be|_{\Spec R \times
U_{\alpha}}$.  Then
\beq
\label{4j}
\be'|_{U_{\alpha}} = \co(D)|_{U_{\alpha}} \oplus \La(D)|_{U_{\alpha}} \otimes
{\cal I}_{(\varepsilon)}.
\eeq
The section $\bp$ is then of the form $\phi + \varepsilon\psi_{\alpha}$ for
some $\psi \in C^0(E)$, and the transition functions are $1 + \varepsilon
f_{\alpha\beta}$ for some $f \in C^1(\End_0 E)$.  The latter hence act as 1 on
the second factor of \re{4j}.

Now decompose $\psi_{\beta} = \psi^{{\cal O}(D)}_{\beta} +
\psi^{\La(-D)}_{\beta}$ and $f_{\alpha\beta} = f^{{\cal O}(D)}_{\alpha\beta} +
f^{\La(-D)}_{\alpha\beta}$ corresponding to the splitting on $U_{\beta}$.  If
$\be'$ is restricted to $(\varepsilon ) \times U_{\alpha}$, then $\varepsilon
\psi^{{\cal O}(D)}_{\beta} = 0$ and $\varepsilon  f^{{\cal O}(D)}_{\alpha\beta}
= 0$, since everything divisible by $\varepsilon $ is now set to zero.
However, $\varepsilon  \psi^{\La(-D)}_{\beta}$ and $\varepsilon
f^{\La(-D)}_{\alpha\beta}$ are not necessarily zero, since not everything in
their images is divisible by $\varepsilon $   {\em in the module} $\La(-D)
\otimes {\cal I}_{(\varepsilon )}$.  Hence $\bp_{(\varepsilon )} = \phi +
\varepsilon  \psi^{{\cal O}(D)}_{\beta}$ on $U_{\beta}$, and $\be_{(\varepsilon
)}$ has transition functions $\left( \begin{array}{cc} 1 &
f^{\La(-D)}_{\alpha\beta} \\ 0 & 1 \end{array} \right)$ with respect to the
splitting \re{4j}.  In other words, the extension class of $E' =
\be_{(\varepsilon )}$ is the projection of $f \in C^1(\End_0 E)$ to
$C^1(\La(2D))$, and the lifting of $\phi'$ is the projection of $\psi \in
C^1(E)$ to $C^1(\La(-D))$.  Hence $(E', \phi')$ is the bundle over the image of
\ephi\ in $\Pj^-_i$ in the family of \re{4h}.  By uniqueness of families
\re{3n} this means that $(\be', \bp')|_{E^+_i \times X}$ is the pullback of the
family of \re{4h}.  \fp

There is a result similar to \re{4l} for the inverse map $M_{i-1} - \Pj W^-_i
\row M_i - \Pj W^+_i$.  Let $\miotm$ be the blow-up of $M_{i-1}$ at $\Pj
W^-_i$.  Hence by \re{4q} the exceptional divisor is $E^-_i = \Pj W^-_i \oplus
\Pj W^+_i$, and $\co_{E^-_i}(E^-_i) = \co(-1,-1)$.  Note that there is an
isomorphism $E^-_i \leftrightarrow E^+_i$.

\begin{propn}
\label{4m}
There is a map $\miotm \row M_i$ such that the following diagram commutes:
$$\begin{array}{ccccc} M_{i-1} - \Pj W^-_i & \lrow & \miotm & \longleftarrow &
E^-_i \\
                  \Big\updownarrow & & \Big\downarrow & & \Big\downarrow \\
                  M_i - \Pj W^+_i & \lrow & M_i & \longleftarrow & \Pj W^+_i.
\end{array}$$

\end{propn}

\pf.  Let $(\be, \bp) \row \miotm \times X$ be the pullback of the universal
family.  We will construct a new family $(\be', \bp')$ of pairs all of which
are in $M_i$.

By uniqueness of families \re{3n}, $(\be, \bp)|_{E^-_i \times X}$ is the
pullback of the family over $\Pj W^-_i$ constructed in \re{4h}.  Thus there is
a surjective sheaf map $\be \row  \co_{E^-_i \times X}(-\Delta)$.  This time
the map does not necessarily annihilate $\bp$.  However, if we tensor by
$\co(E^-_i)$, then the twisted map $\be(E^-_i) \row \co_{E^-_i \times
X}(\Delta)(-1,-1)$ of course annihilates $\bp(E^-_i)$.  If we define $\be'$ to
be the kernel of this twisted map, so that
$$0 \lrow \be' \lrow \be(E^-_i) \lrow \co_{E^+_i \times X}(\Delta)(-1,-1) \lrow
0,$$
then $\be'$ is locally free, and $\bp(E^-_i)$ descends to $\bp' \in H^0(\be')$.
 The remainder of the proof is analogous to that of \re{4l}.  \fp

At last we come to the goal of all the above work.

\begin{propn}
There is a natural isomorphism $\mitp \leftrightarrow \miotm$ such that the
following diagram commutes:
$$\begin{array}{ccccc}M_i - \Pj W^+_i & \lrow & \mitp & \longleftarrow & E^+_i
\\
                  \Big\updownarrow & & \Big\updownarrow & & \Big\updownarrow \\
                  M_{i-1} - \Pj W^-_i & \lrow & \miotm & \longleftarrow &
E^-_i.\end{array}$$

\end{propn}

\pf.  Both $\mitp$ and $\miotm$ are smooth, and by \re{4l} and \re{4m} they
both inject into $M_{i-1} \times M_i$.  Indeed, both injections are embeddings,
since  as is easily checked they annihilate no tangent vectors, and both have
the same image.  This image is precisely the closure of the graph of the
isomorphism $M_i - \Pj W^+_i \leftrightarrow M_{i-1} - \Pj W^-_i$, which proves
the left-hand square; for both $E^-_i$ and $E^+_i$ it is the map $\Pj W^-_i
\oplus \Pj W^+_i \row \Pj W^-_i \times \Pj W^+_i$, which proves the right-hand
square. \fp

{\em Note.}  In light of this result, we will henceforth refer to $\mitp =
\miotm$ simply as $\mitl$, and $E^+_i = E^-_i$ as $E_i$.

Thus $M_i$ is obtained from $M_{i-1}$ by blowing up $\Pj W^-_i$, and then
blowing down the same exceptional divisor in another direction.  Such a blow-up
and blow-down is an example of what is called a {\em flip} in Mori theory.
This paper will not use any of the deep results of Mori theory, but we will see
some of its basic principles in action.

In one case the flip degenerates to an ordinary blow-up.

\begin{propn}
\label{4r}
The moduli space $M_1$ is the blow-up of $M_0 = \Pj H^1(\La^{-1})$ along $X$
embedded via $|K_X\La|$.
\end{propn}

\pf.  Since $W^-_1$ is a line bundle, there is nothing to blow down.  \fp

The other extreme case is also of interest.  Let $w = [(d-1)/2]$, so that $M_w$
is the last moduli space in our sequence.  Let $N$ be the moduli space of
ordinary rank 2 semistable bundles of determinant $\La$.

\begin{propn}
\label{4s}
There is a natural ``Abel-Jacobi'' map $M_w \row N$ with fibre $\Pj H^0(E)$
over a stable bundle $E$.  It is surjective if $d > 2g-2$.
\end{propn}

\pf.  If $i = w$, then $\s \in (0, [d/2] + 1 - d/2)$, so $\s$-stability of
\ephi\ implies ordinary semistability of $E$.  Thus there is a map $M_w \row
N$.  Moreover, ordinary stability of $E$ implies $\s$-stability of \ephi, so
the fibre over a stable $E$ is just $\Pj H^0(E)$.  For $d > 2g-2$, any bundle
$E$ has a nonzero section $\phi$ by Riemann-Roch.  Hence every stable bundle in
$N$ is certainly in the image of $M_w$.  But $M_w$ is complete, so its image is
a complete variety containing the stable set, which must be $N$ itself.  \fp

We may sum up our findings in the following diagram.

\def\arraycolsep{2pt}
$$\begin{array}{ccccccccccccccccc}
 & & \tilde{M}_2 & & & & \tilde{M}_3 & & & & \tilde{M}_4 & & & & \tilde{M}_w &
& \\
 & \swarrow & & \searrow & & \swarrow & & \searrow & & \swarrow & & \searrow &
& \swarrow & & \searrow & \\
M_1 & & & & M_2 & & & & M_3 & & & & \, \cdots \, & & & & M_w \\
\downarrow & & & & & & & & & & & & & & & & \downarrow \\
M_0  & & & & & & & & & & & & & & & &  N
\end{array}$$
\def\arraycolsep{5pt}

All the arrows are birational morphisms except sometimes the one to $N$.

\bit{Their Poincar\'e polynomials}

Before going on to our main application in the next section, let us pause to
see how the flips described above can be used to compute the Poincar\'e
polynomials of our moduli spaces.

\beq
P_t(M_i) = \frac{1}{1-t^2} \Coeff_{x^i} \left(\frac{t^{2d+2g-2-4i}}{xt^4-1} -
\frac{t^{2i+2}}{x-t^2} \right) \left( \frac{(1+xt)^{2g}}{(1-x)(1-xt^2)}
\right).
\eeq

\pf.  Since $\tilde{M}_j$ is the blow-up of $M_{j-1}$ at $\Pj W^-_j$, by the
formula for the Poincar\'e polynomial of a blow-up \cite[p.\ 605]{gh},
$$P_t(\tilde{M}_j) = P_t(M_{j-1}) + P_t(E_j) - P_t(\Pj W^-_j).$$
But $\tilde{M}_j$ is also the blow-up of $M_j$ at $\Pj W^+_j$, so
$$P_t(\tilde{M}_j) = P_t(M_j) + P_t(E_j) - P_t(\Pj W^+_j)$$
as well.  Hence
$$P_t(M_j) - P_t(M_{j-1}) = P_t(\Pj W^+_j) - P_t(\Pj W^-_j).$$
But the Poincar\'e polynomial of any projective bundle splits, so
\beqas
P_t(\Pj W^+_j) - P_t(\Pj W^-_j) & = & P_t(\Pj^{d+g-2-2j})P_t(X_j) -
P_t(\Pj^{i-1})P_t(X_j) \\
& = & \frac{t^{2j}-t^{2d+2g-2-4j}}{1-t^2} P_t(X_j).
\eeqas
A formula for $P_t(X_j)$ was given by Macdonald \cite{mac}:
$$P_t(X_j) = \Coeff_{x^j} \frac{(1+xt)^{2g}}{(1-x)(1-xt^2)}.$$
Hence
$$P_t(M_j) - P_t(M_{j-1}) = \Coeff_{x^j}
\frac{(t^{2j}-t^{2d+2g-2-4j})(1+xt)^{2g}}{(1-t^2)(1-x)(1-xt^2)}.$$
Notice that this formula also produces $P_t(M_0)$ when $j=0$.  So to sum up,
$$\begin{array}{ccl}
P_t(M_i) &  =  & \displaystyle\frac{1}{1-t^2} \Coeff_{x^i} \sum_{j=0}^{i}
\frac{x^{i-j}(t^{2j}-t^{2d+2g-2-4j})(1+xt)^{2g}}{(1-x)(1-xt^2)} \\[15pt]
  & = & \displaystyle\frac{1}{1-t^2} \Coeff_{x^i} \left(
\frac{x^{i+1}-t^{2i+2}}{x-t^2} +
\frac{t^{2d+2g-2-4i}(1-t^{4i-4}x^{i+1})}{xt^4-1} \right)
\left( \frac{(1+xt)^{2g}}{(1-x)(1-xt^2)} \right), \end{array} $$
which agrees with the formula stated after the terms containing $x^{i+1}$ are
removed. \fp

We can use this formula to recover the formula of Harder-Narasimhan \cite{hn}
for the Poincar\'e polynomial of the moduli space $N$ of stable bundles of rank
2, determinant $\La$, and odd degree $d$:

\beq
\label{4u}
P_t(N) = \frac{(1+t^3)^{2g}-t^{2g}(1+t)^{2g}}{(1-t^2)(1-t^4)}.
\eeq

\pf.  When $d > 2g-2$ is odd and $i = w$, then by \re{4s} there is a surjective
map $M_w \row N$ with fibre $\Pj H^0(E)$ over a bundle $E$.  If moreover $d >
4g-4$, then $H^1(E) = 0$ for all stable $E$ (see for example the proof of
\re{3e}), so $M_w$ is then just the $\Pj^{d-2g+1}$-bundle $\Pj (R^0\pi) \be$,
where \be\ is a universal bundle over $N$, and
$$P_t(N) = \frac{1-t^2}{1-t^{2d-4g+4}} P_t(M_w).$$
For simplicity we may as well assume that $d = 4g-3$. Then $w = 2g-2$ and
$$P_t(N) = \frac{1}{1-t^{4g-2}} \Coeff_{x^{2g-2}}\left(\frac{t^{2g}}{xt^4-1} -
\frac{t^{4g-2}}{x-t^2}\right) \left( \frac{(1+xt)^{2g}}{(1-x)(1-xt^2)}
\right).$$
The following argument, due to Don Zagier, then shows that this equals the
Harder-Nara\-sim\-han formula.  Let
$$F(a,b,c,t) = \Coeff_{x^{2g-2}} \frac{(1+xt)^{2g}}{(1-ax)(1-bx)(1-cx)}.$$
Then
$$P_t(N) = \frac{t^{4g-4}F(1,t^2,t^{-2},t) -
t^{2g}F(1,t^2,t^4,t)}{1-t^{4g-2}}.$$
On the other hand,
$$F(a,b,c,t) = \Res_{x=0}\left\{
\frac{x^{1-2g}(1+xt)^{2g}dx}{(1-ax)(1-bx)(1-cx)}\right\};$$
since this has no pole at infinity, by the residue theorem
\beqas
F(a,b,c,t) & = & (-\Res_{x = 1/a} -\Res_{x = 1/b} -\Res_{x = 1/c}) \left\{
\frac{x^{1-2g}(1+xt)^{2g}dx}{(1-ax)(1-bx)(1-cx)}\right\} \\
& = & \frac{(a+t)^{2g}}{(a-b)(a-c)} + \frac{(b+t)^{2g}}{(b-a)(b-c)} +
\frac{(c+t)^{2g}}{(c-a)(c-b)}.
\eeqas
After this substitution, it is a matter of high-school algebra to verify
\re{4u}.  \fp

\bit{Their ample cones}

We now turn to a study of the line bundles over the $M_i$.  Indeed, our goal is
a formula for the dimension of the space of sections of any line bundle over
any $M_i$. Since $M_0$ is just a projective space, the first interesting case
is $M_1$; so we first of all ask what line bundles there are on $M_1$.

\begin{propn}
$\Pic M_1 = \Z \oplus \Z$, generated by the hyperplane $H$ and the exceptional
divisor $E_1$.
\end{propn}

\pf.  Obvious from \re{4r}.  \fp

The case of $M_1$ will be crucial for us, so we introduce the notation
$$\co_1(m,n) = \co((m+n)H - nE_1),$$
$$V_{m,n} = H^0(M_1; \co_1(m,n)).$$  Pushing down to $M_0 = \Pj H^1(\La^{-1})$
then yields $V_{m,n} = H^0(M_0; \co(m+n) \otimes {\cal I}_X^n)$.  That is, an
element of $\Pj V_{m,n}$ is a hypersurface of degree $m+n$ with a singularity
of order $n-1$ at $X$.  The dimension of $V_{m,n}$, which we shall attempt to
calculate, is thus a number canonically associated to $X$, $\La$, $m$, and $n$.

Of course, in many cases this number is easy to compute.  If $m < 0$, for
example, then $V_{m,n} = 0$, since no hypersurface can have a singularity of
order greater than its degree.  If $n < 0$, then $V_{m,n} = H^0(M_0; \co(m+n)
\otimes {\cal I}_X^n) = H^0(M_0; \co(m+n))$, because $\codim X/M_0 = d+g-3 > 1$
by our assumptions on $d$ and $g$, and a section cannot have a pole on a set of
codimension $>1$.  So in this case $\dim V_{m,n} = {m+n+d+g-2 \choose m+n}$.
However, for $m,n \geq 0$, it is quite an interesting problem to calculate
$\dim V_{m,n}$. When $n=1$, these are of course precisely the spaces whose
syzygies are studied by Green and Lazarsfeld \cite{gl}, but for $n > 1$ very
little appears to be known.

What about $M_i$ for $i > 1$?  These give exactly the same information as
$M_1$, for the following simple reason.

\begin{propn}
\label{5a}
For $i > 0$ there is a natural isomorphism $\Pic M_1 \conga \Pic M_i$.
Moreover, if by abuse of notation we denote by $\co_i(m,n)$ the image of
$\co_1(m,n)$ in $\Pic M_i$, then for any $m,n$ there is a natural isomorphism
$V_{m,n} \conga H^0(M_i; \co_i(m,n))$.
\end{propn}

\pf.  By \re{4t}, $M_1$ is isomorphic to $M_i$ except on sets of codimension
$\geq 2$.  Hence divisors, functions, line bundles, and sections can be pulled
back from one to the other and extended over the bad sets in a unique way.  \fp

However, we will certainly not ignore the higher $M_i$ for the rest of the
paper.  Instead, they will be indispensable tools in the study of the
cohomology of $M_1$, to be used as follows.  A naive approach to calculating
$\dim V_{m,n}$ would be to calculate $\x(M_1; \co_1(m,n))$, which is easy using
Riemann-Roch, and then to apply Kodaira vanishing to show that the higher
cohomology all vanished.  This will not work: the hypothesis of Kodaira
vanishing, which is that $K^{-1}_{M_1} \co_1(m,n)$ must be ample, will not
typically be satisfied, and the higher cohomology will not vanish.  But this
problem can be cured by shifting attention to some other $M_i$.  Indeed, under
some mild hypotheses on $m$ and $n$, there will be some $i$ such that
$K^{-1}_{M_i} \co_i(m,n)$ will be ample on $M_i$.  Hence $\dim V_{m,n} =
\x(M_i; \co_i(m,n))$, which will be calculated by an inductive procedure on
$i$.

To carry out this programme, of course, we need to know the ample cone of each
$M_i$.  So our goal in this section will be to prove the following theorem.

\begin{thm}
\label{5b}
For $0 < i < w$, the ample cone of $M_i$ is bounded by $\co_i(1,i-1)$ and
$\co_i(1,i)$.  For $d > 2g-2$, the ample cone of $M_w$ is bounded by
$\co_w(1,w-1)$ and $\co_w(2,d-2)$; for $d \leq 2g-2$, it is bounded on one side
by $\co_w(1,w-1)$, and contains the cone bounded on the other side by
$\co_w(2,d-2)$.
\end{thm}

So as we pass from $i-1$ to $i$, the ample cone flips across the ray of slope
$i-1$, as illustrated for $d = 7$ in the figure.  This is exactly the behaviour
which is predicted by Mori theory; indeed, flips are so named for precisely
this reason.

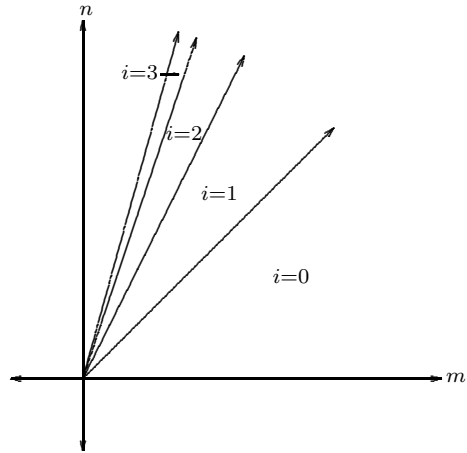
\begin{figure}[hbt]
\begin{center}

\font\thinlinefont=cmr5
\mbox{\beginpicture
\setcoordinatesystem units < 0.500cm, 0.500cm>
\unitlength= 0.500cm
\linethickness=1pt
\setplotsymbol ({\makebox(0,0)[l]{\tencirc\symbol{'160}}})
\setshadesymbol ({\thinlinefont .})
\setlinear
%
%
\linethickness= 0.500pt
\setplotsymbol ({\thinlinefont .})
\putrule from 6.462 17.638 to 6.938 17.638
%
%
\plot 6.811 17.607 6.938 17.638 6.811 17.670 /
%
%
%
\linethickness= 0.500pt
\setplotsymbol ({\thinlinefont .})
\plot 4.398 9.542 6.938 18.749 /
%
%
\plot 6.932 18.488 6.938 18.749 6.809 18.521 /
%
%
%
\linethickness= 0.500pt
\setplotsymbol ({\thinlinefont .})
\plot 4.398 9.542 7.415 18.590 /
%
%
\plot 7.394 18.329 7.415 18.590 7.274 18.370 /
%
%
%
\linethickness= 0.500pt
\setplotsymbol ({\thinlinefont .})
\plot 4.398 9.542 8.684 18.114 /
%
%
\plot 8.628 17.859 8.684 18.114 8.514 17.915 /
%
%
%
\linethickness= 0.500pt
\setplotsymbol ({\thinlinefont .})
\plot 4.398 9.542 11.066 16.209 /
%
%
\plot 10.931 15.985 11.066 16.209 10.841 16.075 /
%
%
%
\linethickness= 0.500pt
\setplotsymbol ({\thinlinefont .})
%
%
\plot 4.335 7.891 4.398 7.637 4.462 7.891 /
\putrule from 4.398 7.637 to 4.398 19.067
%
%
\plot 4.462 18.813 4.398 19.067 4.335 18.813 /
%
%
%
\linethickness= 0.500pt
\setplotsymbol ({\thinlinefont .})
%
%
\plot  2.747 9.605  2.493 9.542  2.747 9.478 /
\putrule from  2.493 9.542 to 13.923 9.542
%
%
\plot 13.669 9.478 13.923 9.542 13.669 9.605 /
%
%
%
\put {$\scriptstyle m$} [lB] at 13.982 9.442
%
%
\put {$\scriptstyle n$} [lB] at 4.239 19.226
%
%
\put {$\scriptstyle i=3$} [lB] at 5.357 17.519
%
%
\put {$\scriptstyle i=2$} [lB] at 6.535 15.892
%
%
\put {$\scriptstyle i=1$} [lB] at 7.473 14.304
%
%
\put {$\scriptstyle i=0$} [lB] at 9.378 12.082
%
%
\linethickness=0pt
\putrectangle corners at  1.873 19.416 and 14.082  7.017
\endpicture}
\end{center}

\end{figure}

The first thing to notice is that, since all the $M_i$ have unique universal
pairs $(\be, \bp) \row M_i \times X$, an expression such as $\det \pi_! \be$,
or $\lii \be_x$ for some $x \in X$, defines line bundles on all the $M_i$,
which agree with one another on the open sets where the maps between different
$M_i$ are defined, and which consequently correspond under the natural
isomorphism of \re{5a}.  Since $\lii \be_x$ and $\det \pi_! \be$ are the
canonical (indeed, essentially the only) examples, we work out what they are on
$M_1$.

\begin{propn}
\label{5d}
On $M_1$, $\lii \be_x = \co_1(0,-1)$ and $\det \pi_! \be = \co_1(-1,g-d)$; that
is, $\co_1(m,n) = \det^{-m} \pi_! \be \otimes (\lii \be_x)^{(d-g)m - n}$.
\end{propn}

\pf.  The universal pair on $M_0 \times X$ is easy to construct directly: it is
the tautological extension
$$0 \lrow \co \lrow \be_0 \lrow \La(-1) \lrow 0$$
determined by the class $id \in \End H^{1}(X;\La^{-1}) = H^{0}(\Pj
H^{1}(\La^{-1}); \co(1)) \otimes H^{1}(X;\La^{-1}) = H^{1}((\Pj H^{1}(\La^{-1})
\times X; \La^{-1}(1))$, together with the constant section $\bp_0 \in
H^0(\co)$.
Recall from \re{4m} that the universal pair $(\be_1, \bp_1) \row M_1 \times X$
is constructed by pulling back $(\be_0, \bp_0)$, twisting by $\co(E^+_1)$, and
modifying at $E^+_1$:
$$0 \lrow \be_1 \lrow \be_0(E^+_1) \lrow \co_{E^+_1 \times X}(\Delta)(-1) \lrow
0.$$
Hence $\lii (\be_1)_x = \lii (\be_0(E^+_1))_x \otimes \co(-E^+_1) = \lii \be_0
\otimes \co(E^+_1) = \co_1(0,-1)$, and
\beqas
\det \pi_! \be_1 & = & \det \pi_! \be_0(E^+_1) \otimes \co((g-2)(E^+_1)) \\
& = & \det \pi_! \co(E^+_1) \otimes \det \pi_! \La(-1)(E) \otimes
\co_1(g-2,2-g) \\
& = & \co_1(1-g,g-1) \otimes \co_1(0,-d-1+g) \otimes \co_1(g-2,2-g) \\
& = & \co_1(-1,g-d).  \fp
\eeqas

The next three results collect some basic information about pullbacks of
$\co_i(m,n)$.

\begin{propn}
\label{5h}
The restriction of $\co_i(m,n)$ to
\begin{tabbing}
{\rm (iii)} \=  \kill
{\rm (i)} \= a fibre of $\Pj W^+_i$ is $\co(n-(i-1)m)$; \\
{\rm (ii)} \= a fibre of $\Pj W^-_i$ is $\co((i-1)m-n)$; \\
{\rm (iii)} \= $f^{-1}(E) \subset M_w$, where $E$ is a stable bundle and $f$ is
the Abel-Jacobi map \\\phantom{xxxxxxxx} of \re{4s}, is $\co(m(d-2)-2n)$.
\end{tabbing}
\end{propn}

\pf.  By \re{4g}, the bundle \be\ in the universal pair restricts to an
extension
$$0 \lrow \co(D) \lrow \be \lrow \La(-D)(-1) \lrow 0$$
on the fibre of $\Pj W^+_i$ over $D \in X_i$.  Hence on this fibre $\lii \be_x
= \co(-1)$ and
$$\det \pi_! \be = \det \pi_! \co(D) \otimes \det \pi_! \La(-D)(-1) =
\co(-\x(\La(-D))) = \co(-d+g-1+i). $$  So  by \re{5d} $\co_i(m,n)$ restricts to
$\co((d-g+1-i)m -(d-g)m + n) = \co((1-i)m +n)$, which proves (i).  Similarly by
\re{4h}, \be\ restricts to an extension
$$0 \lrow \La(-D)(1) \lrow \be \lrow \co(D) \lrow 0$$
on the fibre of $\Pj W^-_i$ over $D \in X_i$.  Hence $\lii \be_x = \co(1)$ and
$$\det \pi_! \be = \det \pi_! \La(-D)(1) \otimes \det \pi_! \co(D) =
\co(\x(\La(-D))) = \co(d-g+1-i). $$
So the previous situation is reversed, and $\co_i(m,n)$ restricts to
$\co((i-1)m-n)$, which proves (ii).  Finally, on a fibre $\Pj H^0(E)$ of the
Abel-Jacobi map, the universal pair restricts to $E(1)$ with the tautological
section.  Hence on this fibre $\lii \be_x = \co(2)$ and $\det \pi_! \be =
\co(d+2-2g)$.  So by \re{5d} $\co_i(m,n)$ restricts to $\co((2g-2-d)m
+2((d-g)m-n)) = \co(m(d-2)-2n)$, which proves (iii).  \fp

\begin{cor}
\label{5j}
On $\mitl$, $\co_i(m,n) = \co_{i-1}(m,n)(((i-1)m-n)E_i)$.
\end{cor}

\pf.  Certainly $\co_i(m,n)$ and $\co_{i-1}(m,n)$ are isomorphic away from
$E_i$, so $\co_i(m,n)$ $= \co_{i-1}(m,n)(q E_i)$ for some $q$.  But
$\co_i(m,n)$ must be trivial on the fibres of $\Pj W^-_i$, and $\co_{E_i}(qE_i)
= \co(-q,-q)$, so by \re{5h}(ii) $q = (i-1)m-n$. \fp

\begin{cor}
\label{5k}
For an effective divisor $D$, let $\iota_D$ be the inclusion of moduli spaces
defined in \re{3p}.  Then $\iota_D^*\co_i(m,n) = \co_i(m,n - m |D|)$.
\end{cor}

\pf.  Choose $x \in X-D$.  Then from \re{5d} and the long exact sequence in
\re{3p}, $\co_i(0,-1) = \lii \be^{\La}_x = \lii (\iota^* \be^{\La(2D)}_x) =
\iota^* \co_i(0,-1)$.  Likewise,
\beqas
\co_i(-1,g-d) & = & \det \pi_! \be^{\La} \\
& = & \det \pi_! \iota^* \be^{\La(2D)} \otimes \det^{-1} \pi_!
\co_D(\be^{\La(2D)}) \\
& = & \det \pi_! \iota^* \be^{\La(2D)} \otimes \bigotimes_{x \in D} (\lii
\be^{\La(2D)}_x)^{-1} \\
& = & \iota^* \co_i(-1, g-d-2|D|) \otimes \iota^* \co_i(0, |D|) \\
& = & \iota^* \co_i(-1, g-d-|D|).  \fp
\eeqas

We now pause to apply these ideas to compute the Picard group of the  moduli
space $N$ of ordinary semistable bundles of determinant $\La$:

\beq
\Pic N = \Z.
\eeq

\pf.  If $g=2$ and $d$ is even, then $N = \Pj^3$ \cite{nr}, so the result is
obvious.  Otherwise, the complement of the stable set $N_s \subset N$ has
codimension $\geq 2$; since $N$ is normal \cite{dn}, this implies $\Pic N_s =
\Pic N$.

By \re{4s} the Abel-Jacobi map $f: M_w \row N$ has fibre $\Pj H^0(E)$ over a
stable bundle $E$.  Tensoring by a line bundle, we may of course assume $d >
4g-4$.  But then $H^1(E) = 0$ (see for example the proof of \re{3e}), so $\dim
\Pj H^0(E) = d+2g-1$ always  and $f$ is locally trivial over $N_s$.  Hence
$\Pic N_s$ is the subgroup of $\Pic M_w$ whose restriction to each projective
fibre of $f$ is trivial.  By \re{5h}(iii) this consists of the bundles
$\co_w(k,k(d/2-1))$ for $k \in \Z$ (where $k$ is even if $d$ is odd). \fp

Denote by $\co(\Theta)$ the $\Q$-Cartier divisor class such that
$f^*\co(\Theta) = \co_w(1,d/2-1)$. Note that this differs slightly from the
normalization in \cite{dn}. The following is then true for any $d$, not just $d
> 4g-4$:

\beq
\label{5e}
f^*\co(\Theta) \conga \co_w(1,d/2-1).
\eeq

\pf.  True by definition if $d > 4g-4$; follows otherwise from \re{5k}, since
$$\iota_D^* \co_w(1,d/2 + |D| -1) = \co_w(1, d/2-1). \fp $$

Now that we know $\Pic N$, we can make the following definition.

\begin{defn}
The {\em Verlinde vector spaces} are
$$Z_k(\La) = H^0(N; \co(k \Theta)),$$
with the convention that $Z_k(\La) = 0$ if $d$ and $k$ are both odd.
\end{defn}

Verlinde's original papers \cite{dv,v} conjectured a striking formula for the
dimensions of these vector spaces, which has since been proved by several
authors. We will give our own proof in \S7; the first step, however, is the
following result, originally due to Bertram \cite{bert}.

\begin{propn}
\label{5f}
For $d > 2g-2$, there is a natural isomorphism $Z_k(\La) \conga
V_{k,k(d/2-1)}$.
\end{propn}

The proof requires the following lemma.

\begin{lemma}
\label{5i}
Let $M$, $N$ be varieties with $N$ normal, and let $f:M \row N$ be a morphism
which is generically a projective bundle.  Then $f_* \co_M = \co_N$.
\end{lemma}

\pf.  This is essentially Stein factorization.  Let $U \subset N$ be the open
set such that $f: f^{-1}(U) \row U$ is a projective bundle.  Then certainly
$f_* \co_{f^{-1}(U)} = \co_U$, so $N' = \SPEC f_* \co_M$ is birational to $N$.
By construction there is a map $f': M \row N'$ such that $f'_* \co_M =
\co_{N'}$.  On the other hand, since $f_* \co_M$ is a coherent sheaf of
$\co_N$-algebras, the birational morphism $N' \row N$ is finite.  But a
birational finite morphism to a normal variety is an isomorphism---this is
essentially Zariski's main theorem; the proof in \cite[III 11.4]{h} goes
through, or see \cite[III.9]{red}. Hence $N' = N$ and $f_* \co_M = \co_N$.  \fp

\pf\ of \re{5f}.  Recall again from \re{4s} that for $d > 2g-2$, the
Abel-Jacobi map $f: M_w \row N$ is surjective with fibre $\Pj H^0(E)$ over a
stable bundle $E$.  If $U \subset N$ is the set of bundles $E$ such that $E$ is
stable and $\dim H^0(E)$ is minimal, then certainly $f: f^{-1}(U) \row U$ is a
projective bundle; for example it is the descent of a trivial projective bundle
over the $\Quot$ scheme.  Moreover, $N$ is always normal \cite{dn}.  So by
\re{5i}, $f_* \co_{M_w} = \co_{N}$.  Hence $f_*f^* \co(k\Theta) =
\co(k\Theta)$, so that
$$f^*: H^0(N;\co(k \Theta)) \row H^0(M_w; \co_w(k, k(d/2-1)))$$
has inverse $f_*$.  \fp

It is worth mentioning, if not proving, a generalization of this result.  Over
the stable set $N_s \subset N$, let $\be \row N_s \times X$ be a universal
bundle, normalized so that $\lii \be |_{N_s \times \{ x \} } = \co$.
(Actually, such a normalization is impossible for $d$ odd, and $\be$ will not
even exist for $d$ even!  However, the obstructions are all in $\Z/2$, and will
cancel in the cases we are considering; for details see \cite{glue}.)  Then let
$U = (R^0 \pi) \be \row N_s$.

\begin{propn}
For $d > 2g-2$, there is a natural isomorphism
$H^0(N_s; S^{m(d-2)-2n}U(m\Theta)) = V_{m,n}$ unless $g=2$ and $d$ is even.
\end{propn}

{\em Sketch of proof.}  The complement of $f^{-1}(N_s) \subset M_w$ has
codimension $\geq 2$ unless $g=2$ and $d$ is even (in which case $N = \Pj^3$
\cite{nr}), so $V_{m,n} = H^0(f^{-1}(N_s); \co(m,n))$.  Also
$(R^0\pi)\co(m,n)|_{N_s} = S^{m(d-2)-2n}U(m\Theta)$, so
$$ H^0(f^{-1}(N_s); \co(m,n)) = H^0(N_s; S^{m(d-2)-2n}U(m\Theta)) $$
as in the proof of \re{5f}.  \fp

Hence seeking a formula for $\dim V_{m,n}$ can be regarded as seeking a
generalization of the Verlinde formula.

At last we return to the determination of the ample cone of $M_i$.  It can of
course be quite difficult to decide whether a given line bundle on a projective
variety is ample.  However, a geometric invariant theory quotient is naturally
endowed with an ample bundle, which is the descent of the ample bundle used in
the linearization.  So we shall work out how the line bundles used in the
linearizations of \S1 descend to $M_i$.  Recall that the linearization was some
power of $\co(\x + 2\s, 4\s) \row \Pj \Hom \times \Pj \C^{\x}$, or more
precisely, its pullback to $\Quot(\La) \times \Pj \C^{\x}$, which by abuse of
notation we still denote $\co(\x + 2\s, 4\s)$.  By further abuse of notation we
refrain from worrying about whether $\x + 2 \s$ and $4\s$ are actually
integers.

\begin{propn}
\label{5c}
The bundle $\co(\x + 2\s, 4\s) \row \Quot(\La) \times \Pj \C^{\x}$ descends to
$\co_i(1,d-1-2\s) \row M_i$.
\end{propn}

\pf. As in \S1, let $U \subset \Quot(\La)$ be the set of quotients $\co^{\x}
\row E \row 0$ of determinant $\La$ such that the induced map $\C^{\x} \row
H^0(E)$ is an isomorphism.  If $\co^{\x} \row \be \row 0$ is the universal
quotient over $U \times X$, then as in \re{3m} there is a universal pair
$(\be(1), \bp) \row U \times \Pj \C^{\x} \times X$ descending to the universal
pair $(\be, \bp)$ on each $M_i$.  Hence $\det \pi_! \be(1) \row U \times \Pj
\C^{\x}$ descends to $\det \pi_! \be = \co_i(-1,g-d) \row M_i$, and for any $x
\in X$, $\lii \be(1)_x \row U \times \Pj \C^{\x}$ descends to $\lii \be_x =
\co_i(0,-1) \row M_i$.

By \cite[III Ex.\ 12.6(b)]{h} $\Pic (U \times \Pj\C^{\x}) = \Pic U \oplus \Pic
\Pj \C^{\x}$.  So to determine a bundle on $U \times \Pj\C^{\x}$, it suffices
to determine it on $\{ E \} \times \Pj \C^{\x}$ and $U \times \{ \phi \}$ for
some $E \in U$, $\phi \in \Pj \C^{\x}$.

On $\{ E \} \times \Pj \C^{\x}$, $\be(1) \conga E(1)$, so $\det \pi_! \be(1) =
\co(\x)$ and $\lii \be_x = \co(2)$.  On $U \times \{ \phi \}$, $\be(1) \conga
\be$, so $\det \pi_! \be(1) = \det \pi_! \be$.  But for all $E \in U$, $H^0(E)
= H^0(\co^{\x})$ and $H^1(E) = 0$.  Consequently $\det \pi_! \be = \co$.
Moreover, there is a canonical map
$$\lii \C^{\x} = \lii H^0(\co(\x)) \lrow \lii H^0(\be) \lrow H^0(\lii \be), $$
so the pullback of $\co(1) \row \Pj \Hom (\lii \C^{\x}, H^0(\La))$ to $U$, also
denoted by $\co(1)$, is precisely $(R^0\pi)\Hom(\La, \lii \be)$.  This is
clearly isomorphic to $\lii \be_x \conga \Hom(\La, \lii \be)_x$, since
$\Hom(\La, \lii \be)$ is trivial on every fibre of $\pi$.

Putting it all together, we find that $\co(0, \x)$ descends to $\co_i(-1,g-d)$
and $\co(1,2)$ descends to $\co_i(0,-1)$.  The result follows after a little
arithmetic.  \fp

{\em Proof} of \re{5b}.  For any $\s \in (\max (0, d/2 - i - 1), d/2 - i)$,
the quotient of $U \times \Pj \C^{\x}$ by the action of $\slx$, linearized by
$\co(\x + 2\s, 4\s)$, gives the same quotient $M_i$.  Hence the descent of
$\co(\x + 2\s, 4\s)$ to $M_i$ is ample for any $\s$ in that interval.  By
\re{5c} and a little arithmetic these bundles span exactly the cones in the
statement of \re{5b}.  Hence those cones are contained in the ample cones of
the $M_i$.  It remains to show that no bundles over $M_i$ outside those cones
are ample, except possibly on one side for $i=w$ and $d \leq 2g-2$.

By \re{5h}(i), the restriction of $\co_i(m,n)$ to a fibre of $\Pj W^+_i$ is
$\co(n-(i-1)m)$.  So $\co_i(m,n)$ can only be ample over $M_i$ if this is
positive, that is, if $(i-1)m < n$.  Thus one side of the ample cone of $M_i$
is where it should be.

Likewise by \re{5h}(ii) the restriction of $\co_{i-1}(m,n) \row M_{i-1}$ to a
fibre of $\Pj W^-_i$ is $\co((i-1)m-n)$.  So for $1 < i \leq w$, that is, when
the dimension of this fibre is positive, $\co_{i-1}(m,n)$ can only be ample
over $M_{i-1}$ if $(i-1)m > n$.  Thus the other side of the ample cone of
$M_{i-1}$ is where it should be.

The only case we have not yet treated is the other side of the ample cone of
$M_w$ for $d > 2g-2$.  In that case there is by \re{4s} a surjective map $M_w
\row N$ onto the moduli space of stable bundles of determinant $\La$.  It is
not an isomorphism, since for example $\Pic M_w = \Z \oplus \Z$ while $\Pic N =
\Z$.  Hence the pullback of the ample bundle $\co(2 \Theta) \row N$ is nef but
not ample, that is, it is in the boundary of the ample cone.  But by \re{5e}
this is precisely $\co(2,d-2)$.  \fp
\bit{Their Euler characteristics}

Now that we know the ample cones of the $M_i$, we can calculate $\dim V_{m,n}$
following the programme outlined in the last section.  We first need a formula
for the canonical bundle of $M_i$:

\beq
\label{6a}
K_{M_i} = \co_i(-3,4-d-g).
\eeq

\pf.  Clearly the canonical bundle is preserved by the isomorphism of \re{5a},
so it suffices to work it out on $M_1$.  But this is easy using \re{4r} and the
standard formulas for the canonical bundle of projective space and of a
blow-up.  \fp

\begin{propn}
\label{6b}
Suppose that $m,n \geq 0$ and that $m(d-2) - 2n > -d+2g-2$.  Let $b =
\left[\frac{n+d+g-4}{m+3}\right] + 1$.  Then $\dim V_{m,n} = \x(M_b;
\co_b(m,n))$.
\end{propn}

The idea of the proof is that $\dim V_{m,n}$ will be an Euler characteristic by
Kodaira vanishing provided that $\co(m,n)$ lies inside some cone in the
translate of the ample fan by $K$.  This is illustrated in the figure for the
case $d=7$.

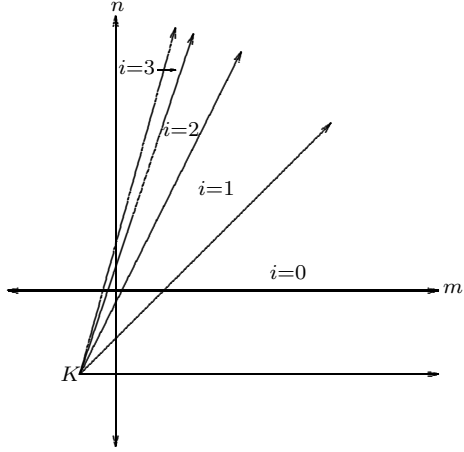
\begin{figure}[htb]
\begin{center}

\font\thinlinefont=cmr5
\mbox{\beginpicture
\setcoordinatesystem units < 0.500cm, 0.500cm>
\unitlength= 0.500cm
\linethickness=1pt
\setplotsymbol ({\makebox(0,0)[l]{\tencirc\symbol{'160}}})
\setshadesymbol ({\thinlinefont .})
\setlinear
%
%
\linethickness= 0.500pt
\setplotsymbol ({\thinlinefont .})
\putrule from 6.462 17.638 to 6.938 17.638
%
%
\plot 6.811 17.607 6.938 17.638 6.811 17.670 /
%
%
%
\linethickness= 0.500pt
\setplotsymbol ({\thinlinefont .})
\plot 4.398 9.542 6.938 18.749 /
%
%
\plot 6.932 18.488 6.938 18.749 6.809 18.521 /
%
%
%
\linethickness= 0.500pt
\setplotsymbol ({\thinlinefont .})
\plot 4.398 9.542 7.415 18.590 /
%
%
\plot 7.394 18.329 7.415 18.590 7.274 18.370 /
%
%
%
\linethickness= 0.500pt
\setplotsymbol ({\thinlinefont .})
\plot 4.398 9.542 8.684 18.114 /
%
%
\plot 8.628 17.859 8.684 18.114 8.514 17.915 /
%
%
%
\linethickness= 0.500pt
\setplotsymbol ({\thinlinefont .})
\plot 4.398 9.542 11.066 16.209 /
%
%
\plot 10.931 15.985 11.066 16.209 10.841 16.075 /
%
%
%
\linethickness= 0.500pt
\setplotsymbol ({\thinlinefont .})
%
%
\plot 5.287 7.891 5.351 7.637 5.414 7.891 /
\putrule from 5.351 7.637 to 5.351 19.067
%
%
\plot 5.414 18.813 5.351 19.067 5.287 18.813 /
%
%
%
\linethickness= 0.500pt
\setplotsymbol ({\thinlinefont .})
%
%
\plot  2.747 11.827  2.493 11.764  2.747 11.700 /
\putrule from  2.493 11.764 to 13.923 11.764
%
%
%
\plot 13.669 11.700 13.923 11.764 13.669 11.827 /
\putrule from  4.398 9.542 to 13.923 9.542
%
%
\plot 13.669 9.478 13.923 9.542 13.669 9.605 /
%
%
%
\put {$\scriptstyle m$} [lB] at 13.982 11.715
%
%
\put {$\scriptstyle n$} [lB] at 5.142 19.226
%
%
\put {$\scriptstyle i=3$} [lB] at 5.357 17.519
%
%
\put {$\scriptstyle i=2$} [lB] at 6.535 15.892
%
%
\put {$\scriptstyle i=1$} [lB] at 7.473 14.304
%
%
\put {$\scriptstyle i=0$} [lB] at 9.378 12.082
%
%
%
%
\put {$\scriptstyle K$} [lB] at 3.810 9.366
\linethickness=0pt
\putrectangle corners at  1.873 19.416 and 14.082  7.017
\endpicture}
\end{center}
\end{figure}

\pf\ of \re{6b}.  Note first that the inequality can be rewritten
$$(d/2-1)(m+3) > n+d+g-4,$$
which guarantees that $b \leq [(d-1)/2]$ and hence that $M_b$ exists.

Suppose that $\frac{n+d+g-4}{m+3}$ is not an integer.  Then $b(m+3) > n+d+g-4 >
(b-1)(m+3)$, so $\co_b(m+3, n+d+g-4)$, which by \re{6a} equals
$K_{M_b}^{-1}\co_b(m,n)$, is in the ample cone of $M_b$ by \re{5b}.  The result
then follows from \re{5a} and Kodaira vanishing.

If $\frac{n+d+g-4}{m+3}$ is an integer, then $\co_{b-1}(m+3, n+d+g-4)$ and
$\co_b(m+3, n+d+g-4)$ are merely nef, so Kodaira vanishing does not apply.
Instead, we move up to $\mbt$.  By \re{5i} the 0th direct image of $\co_{\mbt}$
in the projection $\mbt \row M_b$ is $\co_{M_b}$, and by the theorem on
cohomology and base change \cite[III 12.11]{h} the higher direct images vanish,
so for all $j$, $H^j(\mbt; \co_b(m,n)) = H^j(M_b; \co_b(m,n))$.  By \re{6a} and
the standard formula for the canonical bundle of a blow-up, $K_{\mbt} =
\co_b(-3, 4-d-g)((b-1)E_b)$.  Unfortunately $K_{\mbt}^{-1} \co_b(m,n)$ may not
be ample, so Kodaira vanishing still does not apply.  Instead, we make the
following two claims: first, that $H^j(\mbt; \co_b(m,n)) = H^j(\mbt;
\co_b(m,n)((b-2)E_b))$ for all $j$, and second, that $\co_b(m+3,
n+d+g-4)(-E_b)$ is ample on $\mbt$.  The desired result follows immediately
from these claims, since at last Kodaira vanishing applies to
$\co_b(m,n)((b-2)E_b)$.

To prove the first claim, note that for $0 < k < b$, $H^j(E_b;
\co_b(m,n)(kE_b)) = 0$ for all $j$, since $\co_b(m,n)(kE_b)$ is $\co(-k)$ on
each fibre of $\Pj^{b-1} \row E_b \row \Pj W^+_b$, so that every term in the
Leray spectral sequence vanishes. Hence from the long exact sequence on $\mbt$
of
$$0 \lrow \co_b(m,n)((k-1)E_b) \lrow \co_b(m,n)(kE_b) \lrow \co_b(m,n)
\co_{E_b}(kE_b) \lrow 0,$$
we get isomorphisms $H^j(\mbt; \co_b(m,n)((k-1)E_b)) \conga H^j(\mbt;
\co_b(m,n)(kE_b))$.  The first claim follows by induction.

As for the second claim, note that on $\mbt$, the line bundles
$\co_{b-1}(1,b-2)$, $\co_b(1,b-1)$, and $\co_b(1,b)$ (or $\co_b(2,2b-1)$ if $b
= w$) are all nef, since they are pulled back from nef bundles on $M_{b-1}$ or
$M_b$.  It is easy using \re{5j}, the constraints on $m$ and $n$, and a little
arithmetic to check that $\co_b(m+3, n+d+g-4)(-E_b)$ is in the interior of the
cone generated by these three bundles. \fp

We will have to assume in future that
\beq
\label{6o}
m(d-2) - 2n > -d+2g-2,
\eeq
since otherwise there is no analogue of the last result and
$K^{-1}_{M_i}\co_i(m,n)$ may not be ample for any $i$.  However, for $d \geq
2g$, we still get a complete answer to our problem, for the following reason.

\begin{propn}
\label{6p}
For $d \geq 2g$ and $m(d-2) - 2n < 0$, $V_{m,n} = 0$.
\end{propn}

\pf.  By Riemann-Roch $\deg E \geq 2g$ implies $\dim H^0(E) \geq 2$, so for any
stable bundle $E$, by \re{4s} the fibre $f^{-1}(E)$ of the Abel-Jacobi map is a
projective space of positive dimension.  By \re{5h}(iii), the restriction of
$\co_w(m,n)$ to this is $\co(m(d-2)-2n)$, so any section of $\co_w(m,n)$ must
vanish on $f^{-1}(E)$.  Hence it must vanish on the inverse image $f^{-1}(N_s)$
of the stable subset of $N$.  But this is open, so it must vanish everywhere.
\fp

Let $L_i \row X_i$ be the line bundle defined by $L_i = \det^{-1} \pi_!
\La(-\Delta) \otimes \det^{-1} \pi_! \co(\Delta)$.  Also put $q_i = n-(i-1)m$.

\begin{propn}
\label{6c}
The restriction of $\co_{i-1}(m,n)$ to $\Pj W^-_i$ is $L_i^m(-q_i)$.
\end{propn}

\pf.  Easy from \re{5d} and the description of the universal pair over $\Pj
W^-_i$ in \re{4h}.  \fp

Now let $U_i \row X_i$ be the vector bundle $(W^-_i) \oplus (W^+_i)^*$, and
define numbers
$$N_i = \x(X_i; L_i^m \otimes \La^i W^-_i \otimes S^{q_i-i}U_i), $$
with of course the convention that this is zero when $q_i-i < 0$.  On $M_0$,
which is just projective space, make the additional convention that $\co_0(m,n)
= \co(m+n)$.

\beq
\label{6d}
N_0 = \x(M_0; \co_0(m,n)) = {m+n+d+g-2\choose m+n}.
\eeq

\pf.  Since $X_0$ is just a point and $W^-_0 = 0$, $U_0 = (W^+_0)^*$ is just
the vector space $H^1(\La)^*$.  Hence $S^{m+n}U_0 = H^0(M_0;\co_0(m,n))$ with
our conventions and the result follows.  \fp

\begin{propn}
\label{6e}
Let $0 < i \leq b$, and suppose that $m,n \geq 0$ satisfy \re{6o}.  Then
$$\x(M_i;\co_i(m,n)) - \x(M_{i-1};\co_{i-1}(m,n)) = (-1)^i N_i.$$
\end{propn}

\pf.   By \re{5i} the 0th direct image of $\co_{\mitl}$ in the projection
$\mitl \row M_i$ is $\co_{M_i}$, and by the theorem on cohomology and base
change \cite[III 12.11]{h} the higher direct images vanish, so $\x(\mitl;
\co_i(m,n)) = \x(M_i; \co_i(m,n))$.  Likewise $\x(\mitl; \co_{i-1}(m,n)) =
\x(M_{i-1}; \co_{i-1}(m,n))$, so it suffices to work on $\mitl$.

Suppose first that $q_i \leq 0$, so that $N_i = 0$.  For $0 < j \leq -q_i$,
consider the exact sequence
$$0 \lrow \co_{i-1}(m,n)((j-1)E_i) \lrow \co_{i-1}(m,n)(jE_i) \lrow
\co_{i-1}(m,n) \otimes \co_{E_i}(jE_i) \lrow 0.$$
By \re{6c} the restriction of $\co_{i-1}(m,n)$ to $E_i = \Pj W^-_i \oplus \Pj
W^+_i$ is $L_i^m(-q_i,0)$, and $\co_{E_i}(E_i) = \co(-1,-1)$, so the third term
of the exact sequence becomes $\co(-q_i-j,-j)$ and we get
$$\x(\mitl;\co_{i-1}(m,n)(jE_i)) - \x(\mitl;\co_{i-1}(m,n)((j-1)E_i)) = \x(E_i;
L_i^m(-q_i-j,-j).$$
Summing over $j$ and using \re{5j} yields
$$\x(\mitl;\co_i(m,n)) - \x(\mitl;\co_{i-1}(m,n)) = \sum_{j=1}^{q_i} \x(E_i;
L_i^m(-q_i-j,-j).$$
However, for $0 < i \leq b$ and $m,n,d,g \geq 0$, a little high-school algebra
shows $-q_i < d+g-1-2i$.  Hence for all $j$ in the sum above, $0 < j <
d+g-1-2i$, so every term in the Leray sequence of the fibration $\Pj^{d+g-2-2i}
\row E_i \row \Pj W^-_i$ vanishes.  Hence all terms are zero, as desired.

Now suppose $q_i > 0$.  By an argument similar to the one above,
$$\x(\mitl;\co_i(m,n)) - \x(\mitl;\co_{i-1}(m,n)) = \sum_{j=0}^{q_{i-1}}
\x(E_i; L_i^m(-q_i+j,j).$$
Each term of the right-hand side can be evaluated using the Leray sequence of
the fibration $\Pj^{i-1} \times \Pj^{d+g-2-2i} \row E_i \row X_i$.  Because
$-q_i+j < 0 \leq j$, the only nonzero direct image of $L_i^m(-q_i+j,j)$ is the
$i$th, which is just $L_i^m \otimes \La^i W^-_i \otimes S^{-q_i+j-i}(W^-_i)
\otimes S^j(W^+_i)^*$.  Here the factor of $\La^i W^-_i$ comes from Serre
duality, since the isomorphism $\co(-i) \conga K_{\Ps^{i-1}}$ is not canonical
unless the right-hand side is tensored by such a factor.  Hence
$$\x(E_i; L_i^m(-q_i+j,j)) = (-1)^i \x(X_i; L_i^m \otimes \La^i W^-_i \otimes
S^{-q_i+j-i}(W^-_i) \otimes S^j(W^+_i)^*).$$
Of course the right-hand side is zero if $q_i-j-i < 0$, so the sum need only
run up to $q_i-i$.  The result follows because certainly
$$S^{q_i-i}U_i = \bigoplus_{j+0}^{q_i-i}S^{q_i-j}(W^-_i) \otimes S^j(W^+_i)^*.
\fp $$

\begin{propn}
\label{6f}
For $i > b$, $N_i = 0$.
\end{propn}

\pf.  It suffices to show that if $i > b$, then $q_i-i < 0$, that is,
$(m+n)/(m+1) < i$.  But using $m,n \geq 0$, the definition of $b$, and the
inequality \re{6o}, it is a matter of high-school algebra to check $(m+n)/(m+1)
\leq b$.  \fp

\beq
\dim V_{m,n} = \sum_{i = 0}^{\infty}(-1)^i N_i.
\eeq

\pf.  Put together \re{6b}, \re{6d}, \re{6e}, and \re{6f}.  \fp

Since each $N_i$ can be evaluated using Riemann-Roch on $X_i$, the right-hand
side depends only on $g$, $d$, $m$, and $n$, not on the precise geometry of $X$
and $\La$.  So even before doing the hard work of the next section, we have
found that $\dim V_{m,n}$ depends only on $g$, $d$, $m$, and $n$, which is
rather surprising.

\bit{Don Zagier to the rescue}

All of the results in this section (except \re{6h} and \re{6l}) are due to Don
Zagier and were communicated by him to the author.

In this section we will compute the $N_i$, using the Riemann-Roch theorem and
Macdonald's description \cite{mac} of the cohomology ring of $X_i$.  So we
begin with a review of Macdonald's results.  Let $e_i, \dots , e_g, e'_1, \dots
e'_g \in H^i(X;\Z)$ be generators such that the intersection form is $\sum_j
e_j \otimes e'_j$.  Define classes $\xi, \xi' \in H^1(X_i: \Z)$ and $\eta \in
H^2(X_i; \Z)$ as the K\"unneth components of the divisor $\Delta \subset X_i
\times X$, regarded as belonging to $H^2(X_i \times X; \Z)$:
$$\Delta = \eta + \sum_j (\xi'_j e_j - \xi_j e'_j) + iX. $$
These generate the ring $H^*(X_i;\Z)$.  Moreover, if we put $\s_j = \xi_j
\xi'_j$, then for any multiindex $I$ without repeats,
\beq
\langle \eta^{i-|I|} \s_I , X_i \rangle = 1.
\eeq
This implies that for any two power series $A(x)$, $B(x)$,
\beqa
\langle A(\eta) \exp (B(\eta)\s), X_i \rangle & = & \sum_{k=0}^{\infty} \langle
A(\eta) B(\eta)^k \s^k/k!, X_i \rangle \nonumber \\
& = & \sum_{k=0}^{g} {g \choose k} \Res_{\eta = 0} \left\{ \frac{A(\eta)
B(\eta)^k}{\eta^{i-k+1}} d\eta \right\} \nonumber \\
& = & \Res_{\eta = 0} \left\{ \frac{A(\eta)(1+\eta B(\eta))^g}{\eta^{i+1}}
d\eta \right\} ,
\label{6j}
\eeqa
where $\s = \sum_j \s_j$.  Note that since $\s_j^2 = 0$, $\s^k/k!$ is just the
$k$th symmetric polynomial in the $\s_j$.

Since we will be doing Riemann-Roch, we need to know the Todd class of $X_i$;
luckily this can be worked out in a useful form.

\beq
\label{6g}
\td X_i = \left( \frac{\eta}{1-e^{-\eta}} \right)^{i-g+1} \exp  \left(
\frac{\s}{e^{\eta}-1} - \frac{\s}{\eta} \right).
\eeq

\pf.  Macdonald \cite{mac} shows that the total Chern class of the tangent
bundle of $X_i$ is
$$
c(X_i) = (1+ \eta)^{i-2g+1}\prod_{j=1}^{g}(1+n-\s_i).
$$
Let $h(x) = x/(1-e^{-x})$, so that
$$\td X_i = h(\eta)^{i-2g+1} \prod_{j+1}^{g} h(\eta-\s_j).$$
Expanding $h(\eta-\s_j)$ in a power series around $\eta$ and using $\s_j^2 =
0$, \beqas
\td X_i & = & h(\eta)^{i-g+1} \prod_{j+1}^{g} \left( 1-\s_j
\frac{h'(\eta)}{h(\eta)} \right) \\
 & = & h(\eta)^{i-g+1} \sum_{k=0}^{\infty} (-1)^k \frac{\s^k}{k!} \left(
\frac{h'(\eta)}{h(\eta)} \right) ^k \\
 & = & h(\eta)^{i-g+1} \exp\left( -\s \frac{h'(\eta)}{h(\eta)}\right) ,
\eeqas
which yields the desired formula.  \fp

\begin{lemma}
\label{6h}
For any line bundle $M \row X$ and any $k \in \Z$,
$$\ch \pi_! M(k\Delta) = ((\deg M + ki + 1 - g) - k^2 \s)  e^{k\eta}.$$
\end{lemma}

\pf.  By Grothendieck-Riemann-Roch
\beqas
\ch \pi_! M(k\Delta) & = & \pi_* \ch M(k\Delta) \td X \\
& = & \pi_* \exp((\deg M + ki)X +k \Xi +k \eta) (1+(1-g)X) \\
& = & \pi_* (1 + (\deg M + ki)X) (1+k \Xi -k^2 \s X)  e^{k\eta}  (1+(1-g)X) \\
& = & ((\deg M + ki + 1 - g) - k^2 \s)  e^{k\eta},
\eeqas
where $\Xi = \sum_j (\xi'_j e_j - \xi_j e'_j)$, so that $\Xi^2 = -2\s X$. \fp

\begin{manynotop}
\label{6l}
\begin{tabbing}
{\bf (\theequation )} \= (iii) \= $\ch(\La^i W^-_i)$ \= = \= \kill
{\bf (\theequation )} \>(i) \> \> $\ch(L_i)$ \' = \> $\exp((d-2i)\eta+2\s)$;\+
\\
(ii) \> \> $\ch(\La^i W^-_i)$ \' = \> $\exp((d-3i+1-g)\eta + 3\s)$; \\
(iii) \> \> $\ch(U_i)$ \' = \> $(d-i+1-2g)  e^{-\eta} + (2g-2)  e^{-2\eta} +
\sum_{j=1}^{g}e^{-\eta-\s_i}$.
\end{tabbing}
\end{manynotop}

\pf.
Since  $L_i = \det^{-1} \pi_! \La(-\Delta) \otimes \det^{-1} \pi_!
\co(\Delta)$, by \re{6h}
$$c_1(L_i) = -c_1(\pi_! \La(-\Delta)) -c_1(\pi_! \co(\Delta)) = (d-i+1-g)\eta +
\s + (-i-1+g)\eta + \s = (d-2i)\eta+2\s,$$
which implies (i).
{}From the exact sequence
$$0 \lrow \La(-2\Delta) \lrow \La(-\Delta) \lrow \co_{\Delta}\La(-\Delta) \lrow
0,$$
it follows that $W^-_i = \pi_! \co_{\Delta}\La(-\Delta) = \pi_! \La(-\Delta) -
\pi_! \La(-2\Delta)$ in $K$-theory.  Hence by \re{6h}
$$
\ch W^-_i = ((d-i+1-g)-\s)e^{-\eta} - ((d-2i+1-g)-4\s)e^{-2\eta}.
$$
In particular
$$c_1(\La^i W^-_i) = c_1(W^-_i) = -(d-i+1-g)\eta -\s +2(d-2i+1-g)\eta +4\s =
(d-3i+1-g)\eta + 3\s,$$
which implies (ii).  Again by \re{6h},
$$\ch (W^+_i)^* = \ch \pi_! \La^{-1}(2\Delta) = ((d-2i+g-1)-4\s) e^{-2\eta}.$$
Hence
\beqas
\ch U_i & = & \ch \: (W^-_i) \oplus (W^+_i)^* \\
& = & ((d-i+1-g)-\s)e^{-\eta} + (2g-2)e^{-2\eta} \\
& = & (d-i+1-2g)e^{-\eta} + (2g-2)e^{-2\eta} + \sum_{j=1}^{g} e^{-\eta-\s_i},
\eeqas
which is (iii).  \fp

\begin{propn}
\label{6i}
$\ch (L_i^m \otimes \La^iW^-_i \otimes S^{q_i-i} U_i)$
$$= \Coeff_{t^{q_i-i}}
\left[ e^{ (m(d-2)-2n)\eta } \exp\left((2m+3)\s-\frac{t\s}{e^{-\eta}-t}\right)
\frac{(e^{-\eta}-t)^{-d+i-1+g}}{(1-t)^{2g-2}} \right]. $$
\end{propn}

\pf. The Chern roots of $S^k U_i$ are the sums of $k$ (not necessarily
distinct)  Chern roots of $U_i$, so by \re{6l}(iii)
\beqas
\sum_{k=0}^{\infty} \ch (S^k U_i) t^k & = & \prod_{
\def\arraystretch{.5}\begin{array}{c}\mbox{\rm \small Chern roots} \\ \mbox{
\rm \small $\alpha$ of $U_i$}\end{array}\def\arraystretch{1}
} \frac{1}{1-te^{\alpha}} \\
& = & \left( \frac{1}{1-te^{-\eta}} \right)^{d-i+1-2g} \left(
\frac{1}{1-te^{-2\eta}} \right)^{2g-2} \prod_{j=1}^{g} \left(
\frac{1}{1-te^{-\eta-\s_j}} \right) \\
& = & \frac{(1-te^{-\eta})^{-d+i-1+g}}{(1-te^{-2\eta})^{2g-2}}
\exp \left( \frac{-t\s}{e^{\eta}-t} \right).
\eeqas
Replacing $t$ by $te^{2\eta}$ and taking coefficients of $t^{q_i-i}$ yields
$$\ch (S^{q_i-i} U_i) =  \Coeff_{t^{q_i-i}} \left[ e^{-2(q_i-i)\eta}\frac{(1-t
e^{\eta})^{-d+i-1+g}}{(1-t)^{2g-2}} \exp\left(\frac{-t\s}{e^{-\eta} - t}
\right) \right] .$$
The result then follows using \re{6l}(i) and (ii) and the pleasing identity
$$m(d-2i) + (d-3i+1-g) - 2(q_i-i) = m(d-2)-2n + (d-i+1-g). \fp $$

We are now ready to perform our Riemann-Roch calculation:

\begin{manynotop}
\label{6k}
\def\baselinestretch{2}
\begin{tabbing}
{\bf (\theequation )} \= = \= xxx \= \kill
\> $N_i$ \' = \> $\langle \ch(L_i^m \otimes \La^i W^-_i \otimes S^{q_i-i}U_i)
\td(X_i), X_i \rangle$  \\[7pt]
\> = \> $\displaystyle \Coeff_{t^{q_i-i}}
\Bigg\langle
e^{ (m(d-2)-2n)\eta } \exp\left((2m+3)\s-\frac{t\s}{e^{-\eta}-t}\right)$
\\[7pt]
\> \> \> $\displaystyle \frac{(e^{-\eta}-t)^{-d+i-1+g}}{(1-t)^{2g-2}}
\left(\frac{\eta}{1-e^{-\eta}} \right)^{i-g+1} \exp \left(
\frac{\s}{e^{\eta}-1} - \frac{\s}{\eta}\right),
X_i \Bigg\rangle $ \\[7pt]
\> = \>
$\displaystyle \Coeff_{t^{q_i-i}}
\Res_{\eta = 0}
\Bigg\{  \frac{e^{((d-2)m-2n)\eta}(e^{-\eta}-t)^{-d+i-1+g}}
{(1+t)^{2g-2}(1-e^{-\eta})^{i+1}}$  \\[7pt]
{\bf (\theequation )} \> \> \>
$\displaystyle
\left(e^{-\eta}+\left(2m+3-\frac{t}{e^{-\eta}-t}\right)(1-e^{-\eta})\right)^g
d\eta \Bigg\}$;
\end{tabbing}
\end{manynotop}
\noindent the first equality by Riemann-Roch, the second by \re{6g} and
\re{6i}, and the third by taking
$$A(x) = \left(\frac{x}{1-e^{-x}}\right)^{i-g+1} e^{((d-2)m-2n)x}
\frac{(e^{-x}-t)^{-d+i-1+g}}{(1+t)^{2g-2}} $$
and $$B(x) = 1/(e^x-1) - 1/x + 2m + 3 - t/(e^{-x}-t)$$
in \re{6j}, then combining $g$th powers.

The term in braces is the product of
$\left(\frac{e^{-\eta}-t}{1-e^{-\eta}}\right)^i$ with something independent of
$i$, so make the substitution
$$y = \frac{e^{-\eta}-t}{1-e^{-\eta}},\,\,\,\,\,\,\,\,\,\, e^{-\eta} =
\frac{1+ty}{1+y},\,\,\,\,\,\,\,\,\,\, 1 - e^{-\eta} = \frac{(1-t)y}{1+y},$$
$$ e^{-\eta} - t = \frac{1-t}{1+y}, \,\,\,\,\,\,\,\,\,\,d\eta =
\frac{(1-t)dy}{(1+y)(1+ty)}.$$
Then the residue in \re{6k} becomes
$$\Res_{y=0}\left\{ \frac{a(y)dy}{y^{i+1}} \right\} = \Coeff_{y^i} a(y) $$
for
$$a(y) = \frac{(1+ty)^{q_{d/2}-1}(1+y)^{-q_{d/2}+d-2g+1}}{(1-t)^{d+g-1}}
\Big(1+(2m+3)(1-t)y - ty^2 \Big)^g.$$
Then since $q_i - i = (m+n) - (m+1)i$,
\beqas
\dim V_{m,n} & = & \sum_{i = 0}^{\infty}(-1)^i N_i \\
& = & \sum_{i = 0}^{\infty}(-1)^i \Coeff_{t^{q_i-i}} \Coeff_{y^{i}}  a(y) \\
& = & \Coeff_{t^{m+n}} \left( \sum_{i = 0}^{\infty} (-t^{m+1})^i \Coeff_{y^i}
a(y)  \right) \\
& = & \Coeff_{t^{m+n}}  a(-t^{m+1}).
\eeqas
Thus we obtain the following theorem.  We repeat the definition of $V_{m,n}$
for convenience.

\begin{thm}
\label{6m}
Let $X$ be embedded in $\Pj H^1(\La^{-1})$ via the linear system $|K_X\La|$.
For any $m,n \geq 0$, let $V_{m,n} = H^0(\Pj H^1(\La^{-1}); \co(m+n) \otimes
{\cal I}_X^n)$.  Define
$$F(t) = \frac{(1-t^{m+2})^{-h-1}(1-t^{m+1})^{-h'-1}}{(1-t)^{d+g-1}t^{m+n}}
\Big(1-(2m+3)(1-t)t^{m+1}-t^{2m+3} \Big)^g,$$
where $h = (d-2)m-2n$ and $h' = -h -d +2g -2$.  Then if $m(d-2) - 2n >
-d+2g-2$,
$$\dim V_{m,n} = \Res_{t=0} \left\{ \frac{F(t) dt}{t} \right\} ,$$
that is, the constant term in the Laurent expansion of $F(t)$ at $t=0$.
Moreover, if $d \geq 2g$ and $m(d-2) - 2n < 0$, then $V_{m,n} = 0$.  \fp
\end{thm}

This is the most explicit formula for $\dim V_{m,n}$ we will obtain in general.
 However, in some cases we could obtain completely explicit formulas.  If $m+n$
is small, for example, we could calculate directly, since we would then be
looking at the residue of a function with a pole of low order; for fixed $m+n$,
we would get an explicit polynomial in $g$, $d$, $m$, and $n$.  Otherwise, we
can still use the residue theorem, which says that the sum of the residues at
all the poles of $F(t) dt/t$ is zero.  These poles are of five possible kinds:
$t=0$, $t=\infty$, $t=1$, $t^{m+1}= 1$ but $t \neq 1$, and $t^{m+2}= 1$ but $t
\neq 1$ (note that the last two cases are disjoint).  But in fact $t=1$ is
never a pole, since at that point $1-(2m+3)(1-t)t^{m+1} -t^{2m+3}$ has a triple
zero, and hence the order of $F(t)$ is
$$(-h-1)+(-h'-1)-(d+g-1)+3g = 1 \geq 0.$$
Also, it is straightforward to check that $F(1/t) = -F(t)$, which implies that
$$\Res_{t = \infty}\left\{\frac{F(t)dt}{t}\right\} = \Res_{t =
0}\left\{\frac{F(t)dt}{t}\right\}.$$
Hence
\beq
\label{6q}
-2\dim V_{m,n} = \Bigg( \sum_{\stackrel{\scriptstyle \zeta^{m+1} = 1}{\zeta
\neq 1}} \Res_{t=\zeta} + \sum_{\stackrel{\scriptstyle \zeta^{m+2} = 1}{\zeta
\neq 1}} \Res_{t=\zeta}\Bigg)\left\{\frac{F(t)dt}{t}\right\}.
\eeq
There are poles at the $(m+2)$th roots of unity if and only if $h \geq 0$, and
at the $(m+1)$th roots of unity if and only if $h' \geq 0$.  Thus $\dim
V_{m,n}$ is a sum over the residues at the $(m+2)$th roots if $h' < 0 \leq h$,
a sum over the residues at the $(m+1)$th roots if $h < 0 \leq h'$, and is 0 if
$h, h' < 0$.  (Note that this last case agrees with \re{6p}.)  For $h \geq 0$
it is necessary to calculate the residue of a function with a pole of order
$1+h$, which gets more and more difficult as $h$ grows.  However, when $h=0$,
the calculation is easy, and we can prove the celebrated Verlinde formula.

\beq
\label{6n}
\dim Z_k(\La) = \left( \frac{k+2}{2} \right)^{g-1}
\sum_{j=1}^{k+1}\frac{(-1)^{d(j+1)}}{( \sin \frac{j\pi}{k+2})^{2g-2}}.
\eeq

\pf.  If $d$ and $k$ are both odd, then on symmetry grounds the right-hand side
is zero as desired.  So assume $d$ and $k$ are not both odd.  By \re{5f} $\dim
Z_k(\La) = \dim V_{k,k(d/2-1)}$ for any $d > 2g-2$.  Then $h=0$ and $h' < 0$,
so by \re{6q}
\beqas
\lefteqn{-2 \dim V_{k,k(d/2-1)}} \\ & = & \sum_{\stackrel{\scriptstyle
\zeta^{k+2} = 1}{\zeta \neq 1}} \Res_{t=\zeta}  \left( \frac{-dt/t}{t^{k+2}-1}
\right) \frac{(1-\zeta^{-1})^{d-2g+1}}{(1-\zeta)^{d+g-1}\zeta^{kd/2}}
\Big(1-(2k+3)(\zeta^{-1}-1) -\zeta^{-1}\Big)^g.
\eeqas
But $(1-(2k+3)(\zeta^{-1}-1) -\zeta^{-1}) = (2k+4)(1-\zeta^{-1})$, the residue
is $-1/(k+2)$, and
$$\frac{(1-\zeta^{-1})^d}{(1-\zeta)^d\zeta^{kd/2}} =
\frac{(1-\zeta^{-1})^d}{(1-\zeta)^d\zeta^{-d}\zeta^{(k+2)d/2}} = (-1)^d
\zeta^{(k+2)d/2},$$
so
\beqas
\dim V_{k,k(d/2-1)} & = & (2k+4)^{g-1} \sum_{\stackrel{\scriptstyle \zeta^{k+2}
= 1}{\zeta \neq 1}} (-1)^{d} \zeta^{(k+2)d/2}
\left(\frac{-\zeta}{(1-\zeta)^2}\right)^{g-1} \\
& = & \half (2k+4)^{g-1}
\sum_{\stackrel{\scriptstyle{\xi^{2k+4}=1}}
{\scriptstyle{\xi \neq \pm 1}}}
\frac{(-1)^{d+g-1}\xi^{(k+2)d}}{(\xi^{-1} - \xi)^{2g-2}},
\eeqas
which is equivalent to the Verlinde formula. \fp
\bit{Relation with Bertram's work}

In this appendix we explain briefly, without proving anything, how this paper
is related to Bertram's work on secant varieties.

In \cite{bert}, Bertram considers how to resolve the rational map $\Pj
H^1(\La^{-1}) \row N$.  He shows that blowing up first $X \subset \Pj
H^1(\La^{-1})$, then the proper transform of each of its secant varieties in
turn, produces after $[(d-1)/2]$ steps a smooth variety $\tilde{\Pj}$ having a
morphism to $N$ that agrees with the rational map away from the blow-ups.  The
existence of the morphism is proved by constructing a sequence of families of
bundles, each obtained by an elementary transformation of the last, starting
with the pullback of the tautological family on $\Pj H^1(\La^{-1}) \times X$,
and ending with a family of bundles that are all semistable.  Bertram's
families of bundles can be interpreted, after some twisting, as families of
pairs in our sense, and it follows that his $\tilde{\Pj}$ dominates all of the
$M_i$.  In other words, he performs all of our blow-ups but none of our
blow-downs.  In particular, our blow-up loci are birational to his, that is,
our $\Pj W^-_i$ in $M_{i-1}$ is the proper transform of the $i$th secant
variety in $\Pj H^1(\La^{-1}) = M_0$.  This makes sense, since both are
essentially $\Pj^{i-1}$-bundles over $X_i$.

However, this correspondence is a little more delicate than it seems, because
the $\Pj^{i-1}$-bundles are different: ours is $\Pj W^-_i = \Pj (R^0 \pi)
\co_{\Delta}\La(-\Delta)$, but as Bertram explains, the secant variety is the
image in $\Pj H^1(\La^{-1})$ of $\Pj (R^0 \pi) \co_{\Delta} K\La$.  How is one
projective bundle transformed into another?  If we pull back the lower secant
varieties to $\Pj (R^0 \pi) \co_{\Delta} K\La$ we find that blowing them up and
down induces a {\em Cremona transformation} on each fibre of the projective
bundle.   For example, consider the $\Pj^2$ fibre over $x_1 + x_2 + x_3 \in
X_3$ of the 3rd secant variety.  This of course meets $X \subset \Pj
H^1(\La^{-1})$ in the 3 points $x_1, x_2, x_3$, so if $X$ is blown up, then
$\Pj^2$ gets blown up at those 3 points.  The proper transform of the 2nd
secant variety meets this blown-up $\Pj^2$ in the proper transforms of the 3
lines between the points, so blowing it up does nothing, and blowing it down
blows down the 3 lines.  All in all we have blown up the vertices of a triangle
in the plane, then blown down the proper transforms of the edges.  This is
well-known to recover $\Pj^2$ \cite[V 4.2.3]{h}; indeed it is given in
coordinates by $[z_0, z_1, z_2] \mapsto [z_1 z_2, z_0 z_2, z_0 z_1]$.

If we do the same thing to $\Pj^3$, we find ourselves blowing up the vertices
of a tetrahedron, then blowing up and down---that is to say, flipping---the
proper transforms of the edges, and finally blowing down the proper transforms
of the faces.  Notice that by the time we get to the faces, they have already
undergone Cremona transformations themselves.  More generally, starting with a
simplex in $\Pj^n$, we may flip all of the subsimplices, starting with the
vertices and working our way up.  The varieties we obtain thus fit into a
diagram shaped exactly like that at the end of \S3.  It is not so well-known
that this recovers $\Pj^n$, or that it is given in coordinates by $[z_i]
\mapsto [z_0 \cdots z_{i-1} z_{i+1} \cdots z_n]$, but these facts can be proved
using the theory of toric varieties.

Even that is not quite the end of the story, since over divisors in $X_i$ with
multiple points the transformations are somewhat different.  Over $2x_1 + x_2
\in X_3$, for example, we want to blow up one reduced point and one doubled
point, then blow down one reduced line and one doubled line.  In coordinates,
this is $[z_0,z_1,z_2] \mapsto [z_0^2, z_0 z_1, z_1 z_2]$.  It is an amusing
exercise to work out coordinate expressions for the Cremona transformations
over other divisors with multiple points.


\begin{thebibliography}{99}
\def\baselinestretch{1}
{\small

\bibitem{acgh}{\sc E. Arbarello, M. Cornalba, P.A. Griffiths, {\rm and} J.
Harris,}
{\it Geometry of algebraic curves,} volume I (Springer-Verlag, 1985).

\bibitem{ab}{\sc M.F. Atiyah {\rm and} R. Bott,} The Yang-Mills equations over
Riemann surfaces,
{\it Philos.\ Trans.\ Roy.\ Soc.\ London Ser.\ A} 308 (1982) 523--615.

\bibitem{bert}{\sc A. Bertram,} Moduli of rank-2 vector bundles, theta
divisors, and the geometry of curves in projective space,
{\it J. Diff.\ Geom.\ } 35 (1992) 429-469.

\bibitem{bert2}{\sc A. Bertram,} Stable pairs and stable parabolic pairs,
Harvard preprint.

\bibitem{brad}{\sc S. Bradlow,} Special metrics and stability for holomorphic
bundles with global sections,
{\it J. Diff.\ Geom.\ } 33 (1991) 169-213.

\bibitem{bd}{\sc S. Bradlow {\rm and} G. Daskalopoulos,} Moduli of stable pairs
for holomorphic bundles over Riemann surfaces,
{\it Int.\ J. Math.\ } 2 (1991) 477-513.

\bibitem{dv}{\sc R. Dijkgraaf {\rm and} E. Verlinde,} Modular invariance and
the fusion algebra,
{\it Nucl.\ Phys.\ B} (Proc.\ Suppl.) 5B (1988) 87--97.

\bibitem{d}{\sc S.K. Donaldson,} Instantons in Yang-Mills theory, in
{\it The interface of mathematics and particle physics,} ed.\ D.G. Quillen,
G.B. Segal, and Tsou S.T. (Oxford, 1990).

\bibitem{dn}{\sc J.-M. Drezet {\rm and} M.S. Narasimhan,} Groupe de Picard des
vari\'et\'es de modules de fibr\'es semi-stables sur les courbes alg\'ebriques,
{\it Inv.\ Math.\ } 97 (1989) 53-94.

\bibitem{gp}{\sc O. Garcia-Prada,} Dimensional reduction of stable bundles,
vortices and stable pairs, in preparation.

\bibitem{g1}{\sc D. Gieseker,} On the moduli of vector bundles on an algebraic
surface,
{\it Ann.\ of Math.\ } 106 (1977) 45-60.

\bibitem{gh}{\sc P. Griffiths {\rm and} J. Harris,}
{\it Principles of algebraic geometry} (Wiley, 1978).

\bibitem{grot}{\sc A. Grothendieck,} Technique de descente et th\'eor\`emes
d'existence en g\'eom\'etrie alg\'ebrique, IV: Les sch\'emas de Hilbert,
{\it S\'em.\ Bourbaki} 1960-61, Exp.\ 221; reprinted in
{\it Fondements de la g\'eom\'etrie alg\'ebrique} (Secr\'etariat Math., Paris,
1962).

\bibitem{hn}{\sc G. Harder {\rm and} M.S. Narasimhan,} On the cohomology groups
of moduli spaces of vector bundles over curves,
{\it Math.\ Ann.\ } 212 (1975) 215-248.

\bibitem{h}{\sc R. Hartshorne,}
{\it Algebraic geometry} (Springer-Verlag, 1977).

\bibitem{gl}{\sc R. Lazarsfeld,} A sampling of vector bundle techniques in the
study of linear series, in
{\it Lectures on Riemann surfaces,} ed.\ M. Cornalba, X. Gomez-Mont, and A.
Verjovsky (World Scientific, 1989).

\bibitem{mac}{\sc I.G. Macdonald,} Symmetric products of an algebraic curve,
{\it Topology} 1 (1962) 319-343.

\bibitem{red}{\sc D. Mumford,}
{\it The red book of varieties and schemes,} Lecture notes in mathematics 1358
(Springer-Verlag, 1988).

\bibitem{mf}{\sc D. Mumford {\rm and} J. Fogarty,}
{\it Geometric invariant theory,} second enlarged edition (Springer-Verlag,
1982).

\bibitem{nr}{\sc M.S. Narasimhan {\rm and} S. Ramanan,} Moduli of vector
bundles on a compact Riemann surface, {\it Ann.\ Math.\ } 89 (1969) 1201--1208.

\bibitem{new}{\sc P.E. Newstead,}
{\it Introduction to moduli problems and orbit spaces} (Tata Inst., Bombay,
1978).

\bibitem{oss}{\sc C. Okonek, M. Schneider, {\rm and} H. Spindler,}
{\it Vector bundles on complex projective spaces} (Birkh\"auser, 1980).



\bibitem{glue}{\sc M. Thaddeus,}  A finite-dimensional approach to Verlinde's
factorization principle, preprint.

\bibitem{v}{\sc E. Verlinde,} Fusion rules and modular transformations in 2d
conformal field theory,
{\it Nucl.\ Phys.\ B} 300 (1988) 360--376.

}

\end{thebibliography}
\end{document}